\documentclass[a4paper,10pt
]{amsart}
\usepackage{amsmath, amscd}
\usepackage{amssymb}
\usepackage[alphabetic]{amsrefs}
\usepackage[T1]{fontenc}
\usepackage[all]{xy}
\usepackage{mathtools}
\usepackage{footmisc}
\usepackage{tikz-cd}
\usepackage{enumitem}
\usepackage{relsize} 
\usepackage[bbgreekl]{mathbbol} 
\usepackage{amsfonts}
\usepackage{comment}
\usepackage{stmaryrd}
\usepackage[colorlinks=true, hyperindex, linkcolor=blue, citecolor=, pagebackref=false, urlcolor=blue, linktocpage]{hyperref}
\usepackage[color = blue!20, bordercolor = black, textsize = tiny]{todonotes}
\usepackage{graphicx}
\graphicspath{ {./} }

\DeclareSymbolFontAlphabet{\mathbb}{AMSb} 
\DeclareSymbolFontAlphabet{\mathbbl}{bbold} 

\numberwithin{equation}{section}
\newtheorem{theorem}[subsection]{Theorem} 
\newtheorem{corollary}[subsection]{Corollary}
\newtheorem{lemma}[subsection]{Lemma}
\newtheorem{proposition}[subsection]{Proposition}

\theoremstyle{definition}
\newtheorem{definition}[subsection]{Definition}
\newtheorem{remark}[subsection]{Remark}
\newtheorem{example}[subsection]{Example}

\newtheorem{construction}[subsection]{Construction}
\newtheorem{question}[subsection]{Question}

\newcommand{\arxivlink}[1]{\href{http://arxiv.org/abs/#1}{\texttt{arXiv:#1}}}

\title[]{Associative log rings and their bimodules}

\author{Deik van der Horst}
\email{deik@live.nl}

\author{Tommy Lundemo}
\address{Mathematisch Instituut, Universiteit Utrecht, Budapestlaan 6, 3584 CD Utrecht, Netherlands}
\email{t.lundemo@uu.nl}

\begin{document}
\maketitle

\begin{abstract} We propose a definition of associative log rings. The role of the commutative monoid is replaced by a variant of a crossed module whose source is a monoid, and we use the resulting group action to control the failure of commutativity. This extends the classical notion of a log ring, even in the commutative case. Natural examples include algebras of logarithmic differential operators. The definition lends itself to a canonical construction of a log enveloping algebra and hence of log bimodules, and we prove that the resulting second log Hochschild cohomology group classifies log square-zero extensions. Along the way, we complete the logarithmic analog of Quillen's classification of abelian group objects in augmented algebras: In exact augmented associative log rings, abelian group objects recover log bimodules. Our constructions enjoy animated analogs, and we develop the basic theory of associative animated log rings \emph{en route} to our classification of exact log square-zero extensions.
\end{abstract}

\section{Introduction}

The more flexible notions of \'etaleness and smoothness available in log geometry have made it appealing to construct ``log'' versions of localizing invariants, such as ${\rm THH}$ and ${\rm TC}$ (\cites{Rog09, RSS26}). As a guiding principle, these should interact with log \'etale and smooth morphisms in a manner analogous to their non-log counterparts: For example, Rognes' log ${\rm THH}$ (\cite[Definition 8.11]{Rog09}) satisfies log \'etale descent (\cite[Theorem 1.7]{Lun21}). 

While these log versions often prove useful in concrete situations, a major conceptual defect is that they thus far have avoided descriptions as actual localizing invariants: To a log scheme $X$, there is no well-established category of quasi-coherent sheaves ${\rm QCoh}(X)$ to which one can apply e.g.\ ${\rm THH}$ and recover Rognes' notion.   

Current approaches to overcome this problem make use of Efimov's theory of dualizable categories \cite{Efi24}, see \cite[Conjecture 1.6]{Lun26}, \cite{Sch26}, and \cite{ABBS26}. The currently available ideas suggest that the resulting (dualizable $\infty$-)categories of log ``modules'', while conceptually meaningful, also involve a technical \emph{tour de force} which actively mixes the dualizable machinery with ideas from almost mathematics and mirror symmetry (in the spirit of Vaintrob \cite{Vai17}). 

The level of abstraction involved in these constructions makes the resulting ``non-commutative objects'' associated to e.g.\ a log ring $(A, M)$ bear \emph{a priori} little resemblance to the simple algebraic datum one starts with. This leads us to ask:

\begin{question}\label{q:asslog} Is there a concrete, $1$-categorical notion of ``associative log rings''? If yes, do constructions in classical log geometry that depend upon commutativity admit generalizations to this context?
\end{question}

The purpose of this paper is to argue that the answer to Question \ref{q:asslog} is an emphatic ``yes''. Recall that, in the context of classical log geometry, a pre-log ring $(A, M)$ consists of a commutative ring $A$, a commutative monoid $M$, together with a structure map $\alpha \colon M \to (A, \cdot)$ of commutative monoids. In the associative setting, we argue that the following definition answers Question \ref{q:asslog} affirmatively, in which commutative monoids are replaced by a variant of crossed modules:

\begin{definition}\label{def:prelog} Let $M$ be a monoid with a left action of a group $G$ by monoid automorphisms. A \emph{log crossed module} $\varphi \colon M \to G$ is a monoid map which is equivariant with respect to the conjugation action on $G$ (that is, $\varphi(g.m) = g\varphi(m)g^{-1}$) and satisfies the \emph{Peiffer identity} $mm' = (\varphi(m).m')m$. 
    
    A \emph{pre-log ring} $(A, M) = (A, \varphi \colon M \to G, \alpha)$ consists of 
\begin{enumerate}
    \item[-] a log crossed module $\varphi \colon M \to G$;
    \item[-] an associative ring $A$ with a left action of $G$ by ring automorphisms; and
    \item[-] an equivariant monoid map $\alpha \colon M \to (A, \cdot)$ 
\end{enumerate}

\noindent subject to the Peiffer identity $\alpha(m)a = (\varphi(m).a)\alpha(m)$. 
\end{definition}

See Remark \ref{rem:logcross} for some context for the notion of a log crossed module. While we refer to Section \ref{subsec:exprelog} for a series of examples, let us stress already now that this definition extends that of a classically commutative pre-log ring already in the case where everything in sight (i.e.\ $A$, $M$, and $G$) is commutative. For this reason, we shall consistently write ``classically commutative'' for the (pre-)log rings considered in the literature so far, and write ${\rm CCPreLog}$ and ${\rm CCLog}$ for their categories. 

The Peiffer identities give us sufficient control over the lack of commutativity to carry many constructions in classical log geometry that seem to depend upon commutativity over to this associative setting. As a motivating example, this leads to a general notion of log bimodules that fit harmoniously as coefficients for Rognes' log Hochschild theory 
\cite{Rog09}, as we now explain:

\subsection{Log Hochschild homology and log bimodules} If $A$ is an associative ring which is flat (over ${\Bbb Z}$) and $J$ is an $A$-bimodule, it is classical that its Hochschild homology ${\rm HH}_*(A ; J)$, defined as the homology of the cyclic bar construction, is computed by the groups ${\rm Tor}_*^{A^e}(A, J)$, where we recall that left modules over the enveloping algebra $A^e := A \otimes A^{\rm op}$ are equivalent to left $(A, A)$-bimodules. There is an analogous observation for Rognes' log Hochschild homology \cite[Definition 3.23]{Rog09}: If $(A, M)$ is a classically commutative pre-log ring for which $A$ is flat over ${\Bbb Z}[M]$, then there is a natural isomorphism \begin{equation}\label{eq:loghhtor}{\rm HH}_*(A, M) \cong {\rm Tor}_*^{(A, M)^e}(A, A),\end{equation} see e.g.\ \cite[Section 13]{Rog09} or \cite[Proposition 1.4]{BLPO23a}. Here, $(A, M)^e$ is the pushout of the diagram \begin{equation}\label{eq:reppush}A \otimes A \xleftarrow{\alpha \otimes \alpha} {\Bbb Z}[M] \otimes {\Bbb Z}[M] \xrightarrow{{\rm sh}} {\Bbb Z}[M^{\rm gp}] \otimes {\Bbb Z}[M]\end{equation} of commutative rings, where ${\rm sh}(m \otimes n) := \gamma(m) \otimes mn$, with $\gamma \colon M \to M^{\rm gp}$ the canonical group completion map. Observe that ${\rm sh}$ is a ring map precisely when $M$ is commutative. The flatness hypothesis may, of course, be removed once everything in sight is derived. 

This suggests that, for a classically commutative pre-log ring $(A, M)$, the category of \emph{bi}modules could simply be taken to be left modules over $(A, M)^e$. 

\subsection{Log Hochschild cohomology and exact square-zero extensions} Spelling out the definitions, for classically commutative $(A, M)$, a left $(A, M)^e$-module $J$ consists of an $(A, A)$-bimodule with a left $M^{\rm gp}$-action, subject to the compatibility criteria $\alpha(m)j\alpha(n) = \gamma(m)j\alpha(mn)$ and $g(ajb) = a(gj)b$. As far as the authors are aware, only the case of the symmetric bimodule $J = A$ with trivial $M^{\rm gp}$-action has been considered in the literature up to this point. We believe that this has obscured some structure that is very natural to consider in this context, and studying properties of $(A, M)^e$-modules eventually led us to propose Definition \ref{def:prelog}. 

For example, while the definition of $(A, M)^e$ as the pushout of \eqref{eq:reppush} only makes sense in a commutative context, we explain in Section \ref{sec:logenv} that the Peiffer conditions allow us to construct for any pre-log ring $(A, M)$ a log enveloping algebra $(A, M)^e$ which recovers the pushout of \eqref{eq:reppush} in the classically commutative case. In addition, we construct an associative generalization of the \emph{logification} construction (or ``passing to the associated log structure''), which, in the classically commutative case, involves a pushout in \emph{commutative} monoids. Since logification is the main ingredient in the definition of a (strict) log square-zero extension, we may now mimic it in this general associative context. In analogy with \eqref{eq:loghhtor}, we take the groups ${\rm Ext}_{(A, M)^e}^*(A, J)$ as the definition of log Hochschild cohomology (in the spirit of \cites{HHL26a, HHL26b, GHH26}) once $(A, M)^e$ is formed in a derived manner. 

\begin{theorem}\label{thm:loghch} Let $(A, \varphi \colon M \to G)$ be a log ring and let $J$ be a left $(A, M)^e$-module. The group ${\rm Ext}^2_{(A, M)^e}(A, J)$ identifies canonically with that of isomorphism classes of exact square-zero extensions $(B, \psi \colon N \to H) \to (A, \varphi \colon M \to G)$ by $J$. 
\end{theorem}

Here, being \emph{exact} refers to the property of the canonical map $N \to M \times_G H$ being an isomorphism of monoids. This notion also plays an important role in the results leading up to Theorem \ref{thm:loghch}, in which we pin down the infinitesimal theory of associative log rings in the spirit of Quillen \cite{Qui70}: 

\subsection{Abelian group objects in augmented log rings} An often-exploited theme in both algebra and geometry is that infinitesimal changes in multiplicative structures are governed by purely additive structures. Quillen \cite{Qui70} provides one way of making this slogan precise: The abelianization functor \begin{equation}\label{eq:quillenab}{\rm CRing}_{/A} \xrightarrow{{\rm Ab}_+} {\rm Ab}({\rm CRing}_{/A}) \simeq {\rm LMod}(A)\end{equation} from commutative rings augmented over $A$ may be described as ${\rm Ab}_+(B) = A \otimes_B \Omega^1_B$ under the given identification with $A$-modules, with $\Omega^1_B$ the module of differentials. In the associative context, we instead obtain bimodules: The analogous result reads \[{\rm Ring}_{/A} \xrightarrow{{\rm Ab}_+} {\rm Ab}({\rm Ring}_{/A}) \simeq {\rm LMod}(A^e), \quad {\rm Ab}_+(B) \cong A^e \otimes_{B^e} I_{B},\] where $I_{B}$ is the kernel of the map $B^e \to B, b \otimes b' \mapsto bb'$. 

In the classically commutative case, one might hope that some variant of the category ${\rm Ab}({\rm CCLog}_{/(A, M)})$ may be taken as the \emph{definition} of log modules over $(A, M)$, in light of Quillen's results. As many constructions in log geometry pass through an \emph{exactification} procedure, one very natural guess is to consider the category ${\rm Ab}({\rm CCLog}_{/(A, M)}^{\rm ex})$ of abelian group objects whose structure map $(B, N) \to (A, M)$ is exact (cf.\ Rognes--Sagave \cite[Remark 8.8]{Rog09}). 

This suggestion was picked up by the second-named author (cf.\ \cite[Chapter 5]{Lun22} or \cite[Section 2]{Lun25}), whose findings can be compactly summarized as follows: 

\begin{enumerate}
    \item The exactified category ${\rm Ab}({\rm CCLog}^{\rm ex}_{/(A, M)})$ does not depend upon the log structure: Indeed, it is equivalent to ${\rm LMod}(A)$; but
    \item it nonetheless \emph{is} the correct category of modules from a deformation-theoretic point of view: One can formulate a suitable ``exactified abelianization'' \[{\rm Ab}_+^{\rm ex} \colon {\rm CCLog}_{/(A, M)} \to {\rm Ab}({\rm CCLog}^{\rm ex}_{/(A, M)}) \simeq {\rm LMod}(A)\] with the property that ${\rm Ab}_+^{\rm ex}(B, N) \simeq A \otimes_B \Omega^1_{(B, N)}$.
\end{enumerate}

\noindent The following reinforces point (2) by generalizing it to the associative setting: 

\begin{theorem}\label{thm:abgrpob} Let $(A, M)$ be an associative log ring. There is an equivalence of categories ${\rm Ab}({\rm Log}^{\rm ex}_{/(A, M)}) \simeq {\rm LMod}((A, M)^e)$, and an exactified abelianization \[{\rm Ab}_+^{\rm ex} \colon {\rm Log}_{/(A, M)} \to {\rm Ab}({\rm Log}_{/(A, M)}^{\rm ex}) \simeq {\rm LMod}((A, M)^e)\] determined by $(B, N) \mapsto (A, M)^e \otimes_{(B, N)^e} I_{(B, N)}$, where $I_{(B, N)}$ is the kernel of an augmentation map $(B, N)^e \to B$.  
\end{theorem}

As opposed to the classically commutative case, the category of abelian group objects now depends upon the log structure. We consider this as strong evidence that the category ${\rm LMod}((A, M)^e)$ is indeed the correct notion of bimodules over the log ring $(A, M)$.

\subsection{Animated log rings and square-zero extensions} In order for the group ${\rm Ext}^2_{(A, M)^e}(A, J)$ considered in Theorem \ref{thm:loghch} to capture \emph{all} square-zero extensions, as opposed to only those with an additive section, it is necessary to interpret the log enveloping algebra $(A, M)^e$ in a fully derived manner, which we do for the remainder of the introduction. This is analogous to the distinction between $A \otimes A^{\rm op}$ and $A \otimes^{\Bbb L} A^{\rm op}$ in the classical non-logarithmic setting. To this end, we develop the theory of animated (associative) log rings, and make all relevant constructions homotopically meaningful within that category. Analogously to \cite[Theorem 4.17]{SSV16} in the commutative case, the resulting log square-zero extensions harmonize well with the derived algebro-geometric notion: 

\begin{theorem}\label{thm:pullback} Let $(A, M)$ be a log ring and $J$ an $(A, M)$-bimodule. There is an equivalence of anima\[{\rm LogDer}((A, M), J[1]) \xrightarrow{\simeq} {\rm SqZero}^{\rm ex}((A, M), J)\] induced by associating to a log derivation $(d, d^\flat)$ the exact square-zero extension $(B, N) \to (A, M)$ stemming from the pullback square \[\begin{tikzcd}(B, N) \ar{r} \ar{d} & (A, M) \ar{d}{(d, d^\flat)} \\ (A, M) \ar{r}{(d_0, d_0^\flat)} & (A, M) \ltimes J[1]\end{tikzcd}\] of animated log rings, where $(d_0, d_0^\flat)$ is the trivial log derivation. 
\end{theorem}

In the classically commutative case, the above is a very slight strengthening of \cite[Theorem 4.17]{SSV16} in that we do not have integrality hypotheses, though this generalization is also easily obtained using the results of \cite[Section 4]{SSV16} relating exact and strict morphisms in the integral case. Applying this result in the special case of a square-zero extension of discrete log rings, Theorem \ref{thm:loghch} follows readily from Theorem \ref{thm:pullback}, as explained in Section \ref{sec:lastsection}.

\subsection{Acknowledgments} The authors would like to thank Federico Binda, M\'arton Hablicsek, Lennart Meier, Doosung Park, and Steffen Sagave for discussions and comments. This work is partially supported by the NWO-grant VI.Veni.242.129.

\subsection{Use of generative AI} A commercial large language model (LLM) was used throughout the project for advanced literature searches and proofreading, during which the LLM correctly pointed out that the present unconditional version of Proposition \ref{prop:logenvlogi} holds true. We only expected its statement to hold under integrality hypotheses, for reasons we elaborate upon in Remark \ref{rem:logenvlogi}. The presented proof of Proposition \ref{prop:logenvlogi} is due to the authors and differs from that proposed by the LLM. The paper is written by the authors and we assume responsibility for its content. 

\subsection{Outline} We discuss our definition of associative pre-log and log rings in Sections \ref{sec:assprelog} and \ref{sec:asslog}. In Section \ref{sec:logenv}, we construct the log enveloping algebra $(A, M)^e$ for any pre-log ring $(A, M)$. We compare their left modules to an intrinsically defined notion of log bimodules and discuss invariance under passing to the associated log structure (``logification''). In Section \ref{sec:properties} we discuss associative analogs of some often-used properties in classical log geometry. Exactness in particular is used in Section \ref{sec:abgrp} to relate abelian group objects to log bimodules, proving Theorem \ref{thm:abgrpob}. Log Hochschild cohomology, together with its relationship to associative log derivations, is introduced in Section \ref{sec:loghch}, while Theorem \ref{thm:loghch} and Theorem \ref{thm:pullback} are proved in Section \ref{sec:lastsection}. 

\subsection{Conventions} All rings and monoids are assumed unital. Given a group $G$ acting on the left on a monoid $M$ by monoid automorphisms, we shall write $G \ltimes M$ for their semidirect product, with multiplication $(g, m)(g', m') = (gg', m(g.m'))$. While this leads to slightly awkward formulas of the form $(g, m) = (1, m)(g, 1)$, the standard alternatives would eventually lead us to write e.g.\ $J \oplus A$ for split square-zero extensions of a ring $A$, which we consider to be a higher price to pay. 

To ease notation we carry out our constructions over an absolute base (that is, ${\Bbb Z}$ with its trivial pre-log structure). We expect little difficulty in establishing the results herein for a more general pre-log base; for example, Theorem \ref{thm:abgrpob} holds verbatim in light of the equivalence ${\rm Ab}({\rm Log}^{\rm ex}_{/(A, M)}) \simeq {\rm Ab}(({\rm Log}_{(R, P)/})^{\rm ex}_{/(A, M)})$ for any map $(R, P) \to (A, M)$ of pre-log rings.

\section{Associative pre-log rings}\label{sec:assprelog} To begin, we give the necessary background for Definition \ref{def:prelog}, and we develop the basic theory for the resulting category of associative pre-log rings, emphasizing new examples and analogies with the classically commutative case. 

\subsection{Classically commutative pre-log rings} Recall that a (classically commutative) pre-log ring $(A, M)$ consists of a commutative ring $A$, a commutative monoid $M$, and a monoid map $\alpha \colon M \to (A, \cdot)$. 

\subsection{Log crossed modules} In the associative setting, the following shall play the role of the commutative monoid in the definition of a pre-log ring; or rather, the role of the commutative monoid together with its group completion map $\gamma \colon M \to M^{\rm gp}$:  

\begin{definition}\label{def:logcrmod} A \emph{log crossed module} $\varphi \colon M \to G$ consists of a group $G$, a monoid $M$ with a left $G$-action by monoid automorphisms, and a monoid map $\varphi \colon M \to G$ subject to the equivariance condition $\varphi(g.m) = g\varphi(m)g^{-1}$ and the \emph{Peiffer \cite{Pei49} identity} $mm' = (\varphi(m).m')m$. We shall write ${\rm LogCrMod}$ for the resulting category of log crossed modules.
\end{definition}

The above definition is intended to interpolate between the classical notion of a crossed module and that of a crossed semimodule; we discuss the historical context of Definition \ref{def:logcrmod} below:

\begin{remark}\label{rem:logcross} There is a ``monoid''-variant of crossed modules, called \emph{crossed semimodules}. As far as the authors are aware, this notion first appeared in Joyal--Street \cite[Section 2, Example 4]{JS86} and was later studied by A.\ Patchkoria \cite{Pat98}. By definition, a crossed semimodule consists of a monoid map $\varphi \colon M \to N$ with an action of $N$ on $M$ by monoid endomorphisms subject to the conditions $n\varphi(m) = \varphi(n.m)n$ and the Peiffer identity $mm' = (\varphi(m).m')m$. Our Definition \ref{def:logcrmod} is thus exactly a crossed semimodule whose target is a group.

In fact, the exact content of Definition \ref{def:logcrmod} already appears in the literature. Indeed, this structure appears under the name of a ``$G$-crossed monoid'' as a decategorification of the $G$-crossed braided monoidal categories discussed in \cite[Example 1.6]{JPR22}. Since it will be essential for us to consider the category of such crossed monoids where both $M$ and $G$ vary, we decided to give this structure a new name to avoid confusion. 
\end{remark}

\begin{example}[The adjoint log crossed module]\label{ex:ada} Let $A$ be a ring and consider \[{\rm Ad}(A) := \{(a, \sigma) \in A \times {\rm Aut}(A) \mid aa' = \sigma(a')a \text{ for all }  a' \in A\}\] as a monoid with multiplication $(a, \sigma)(b, \gamma) = (ab, \sigma \circ \gamma)$. We consider the left action of ${\rm Aut}(A)$ on ${\rm Ad}(A)$ given by $\gamma.(a, \sigma) := (\gamma(a), \gamma \circ \sigma \circ \gamma^{-1})$. With these definitions, the map \[q \colon {\rm Ad}(A) \to {\rm Aut}(A), \quad (a, \sigma) \mapsto \sigma\] is readily verified to determine a log crossed module. The notation is meant to signify the fact that, once restricted to those pairs $(a, \sigma)$ with $a \in A^\times$, the defining relation reads $\sigma(a') = aa'a^{-1}$ for all $a' \in A$; that is, $\sigma = {\rm Ad}_a$.
\end{example}

\begin{example}[The unit log crossed module] We consider the submonoid of ${\rm Ad}(A)$ consisting of those pairs $(a, \sigma)$ with $a \in A^\times$. Since this forces $\sigma = {\rm Ad}_a$, this submonoid is canonically isomorphic to $A^\times$. The map $A^\times \to {\rm InnAut}(A), a \mapsto {\rm Ad}_a$ is a log crossed module, and is naturally a subobject of the adjoint log crossed module ${\rm Ad}(A) \to {\rm Aut}(A)$. 
\end{example}

\begin{example}[Commutative monoids] Let $M$ be a commutative monoid. Then the group completion map $\gamma \colon M \to M^{\rm gp}$ is a log crossed module, with respect to the trivial $M^{\rm gp}$-action on $M$. 
\end{example}

\begin{example}[Permutation symmetries]\label{ex:permsymm} The inclusion ${\Bbb N}^{\oplus n} \to \Sigma_n \ltimes {\Bbb Z}^{\oplus n}$ is a log crossed module, where $\Sigma_n$ acts on ${\Bbb N}^{\oplus n}$ by permuting coordinates. 
\end{example}

\subsection{Pre-log rings} With the notion of a log crossed module in place, we recall the definition of an associative pre-log ring: 

\begin{definition}\label{def:prelogmaintext} A \emph{pre-log ring} $(A, M)$ = $(A, \varphi \colon M \to G, \alpha)$ consists of 
    
    \begin{enumerate} \item an associative ring $A$; \item a log crossed module $\varphi \colon M \to G$; \item an action of $G$ on $A$ by ring automorphisms; and \item an equivariant monoid map $\alpha \colon M \to (A, \cdot)$ 
    \end{enumerate}
\noindent satisfying the Peiffer condition $\alpha(m)a = (\varphi(m).a)\alpha(m)$.
\end{definition}

When the remaining structure is clear from context, we shall sometimes simply write $(A, \varphi \colon M \to G)$ or even $(A, M)$. The role of the group completion map $\gamma \colon M \to M^{\rm gp}$ in the classically commutative situation is now played by the map $\varphi \colon M \to G$. For example, we say that $(A, M)$ is \emph{integral} if $\varphi$ is injective. 

\subsection{The category ${\rm RingRep}$} It is occasionally useful to isolate the resulting category of pairs $(A, G)$ involved in the definition of a pre-log ring. We define a category ${\rm RingRep}$ whose objects are pairs $(A, G)$ where $A$ is an associative ring and $G$ is a group acting on $A$ by ring automorphisms. Morphisms $(f, f^\#) \colon (A, G) \to (B, H)$ in ${\rm RingRep}$ consist of a ring map $f \colon A \to B$, a group map $f^\# \colon G \to H$, subject to the compatibility criterion $f(g.a) = f^\#(g).f(a)$. 

\subsection{Morphisms in ${\rm PreLog}$} In the same way that a classically commutative pre-log ring can be described as a morphism $\alpha \colon M \to (A, \cdot)$ of commutative monoids, one can describe a general associative pre-log ring as a morphism of log crossed modules \[\begin{tikzcd}M \ar{rrrr}{m \mapsto(\alpha(m), a \mapsto \varphi(m).a)} \ar{d}{\varphi} & & & & {\rm Ad}(A) \ar{d}{q} \\ G  \ar{rrrr}{g \mapsto (a \mapsto g.a)} & & & & {\rm Aut}(A).\end{tikzcd}\] ${\rm Hom}$-sets in ${\rm PreLog}$ are defined by the pullback squares \[\begin{tikzcd}[column sep = small]{\rm Hom}_{\rm PreLog}((A, M \to G), (B, N \to H)) \ar{r} \ar{d} & {\rm Hom}_{{\rm RingRep}}((A, G), (B, H)) \ar{d} \\ {\rm Hom}_{\rm LogCrMod}(M \to G, N \to H) \ar{r} & {\rm Hom}_{{\rm Grp}}(G, H) \times {\rm Hom}_{\rm Mon}(M, (B, \cdot)).\end{tikzcd}\] A morphism of pre-log rings is thus a triple \[(f, f^\flat, f^\#) \colon (A, \varphi \colon M \to G) \to (B, \psi \colon N \to H)\] consisting of a ring map $f \colon A \to B$, monoid map $f^\flat \colon M \to N$, and group map $f^\# \colon G \to H$, respecting all involved structure. 

\subsection{Examples of pre-log rings}\label{subsec:exprelog} We record here a series of natural examples: 

\begin{example}[Classically commutative pre-log rings] The category of classical commutative pre-log rings embeds into the category of pre-log rings by sending $(A, M, \alpha)$ to $(A, \gamma \colon M \to M^{\rm gp}, \alpha)$, where $\gamma$ is the canonical map to the group completion $M^{\rm gp}$, assumed to act trivially on $M$ and $A$.
\end{example}

\begin{example}[The adjoint pre-log ring]\label{ex:adaprelog} Let $A$ be an associative ring. With notation as in Example \ref{ex:ada}, $(A, q \colon {\rm Ad}(A) \to {\rm Aut}(A))$ is a pre-log ring with respect to the structure map $\alpha(a, \sigma) = a$. 
\end{example}

\begin{example}[The pre-log ring associated to a log crossed module]\label{ex:zm} If $M$ is a commutative monoid, one often considers the classically commutative pre-log ring $({\Bbb Z}[M], M)$. If $\varphi \colon M \to G$ is a log crossed module, the construction generalizes to the associative pre-log ring $({\Bbb Z}[M], \varphi \colon M \to G)$ with the action induced by that of $G$ on $M$. 
\end{example}

\begin{remark}[The building blocks of a pre-log ring] We find it worthwhile to pause here to record an instructive analogy with the classically commutative context. If $(A, M)$ is a classically commutative pre-log ring, one often exploits the fact that it is the pushout of the diagram \[(A, \{1\}) \xleftarrow{} ({\Bbb Z}[M], \{1\}) \xrightarrow{} ({\Bbb Z}[M], M)\] in the category of classically commutative pre-log rings. This observation uses in particular the fully faithful embedding ${\rm CRing} \subset {\rm CCPreLog}, A \mapsto (A, \{1\})$.  

The analogous embedding in the general associative context should remember the action of $G$ on $A$; that is, \begin{equation}\label{eq:repemb}{\rm RingRep} \subset {\rm PreLog}, \quad (A, G) \mapsto (A, \{1\} \to G),\end{equation} where the pre-log ring $(A, \{1\} \to G)$ has the given action of $G$ on $A$. It is then the case that a pre-log ring $(A, \varphi \colon M \to G)$ may be realized as the pushout of  \[(A, \{1\} \to G) \xleftarrow{} ({\Bbb Z}[M], \{1\} \to G) \xrightarrow{} ({\Bbb Z}[M], \varphi \colon M \to G)\] in ${\rm PreLog}$, where the right-hand term is Example \ref{ex:zm} applied to $\varphi \colon M \to G$, the middle term is \eqref{eq:repemb} applied to its underlying object $({\Bbb Z}[M], G)$ in ${\rm RingRep}$, while the left-hand term is \eqref{eq:repemb} applied to $(A, G)$.  
\end{remark}

\begin{example}[Permutation symmetries of ${\Bbb A}^n$] Building upon Example \ref{ex:permsymm}, we obtain the pre-log ring $({\Bbb Z}[x_1, \dots, x_n], {\Bbb N}^{\oplus n} \to \Sigma_n \ltimes {\Bbb Z}^{\oplus n})$, where $\alpha(e_i) = x_i$. 
\end{example}

\begin{example}[A commutative pre-log ring which is not classical]\label{ex:comm} The triple $({\Bbb Z}[x, y], \varphi \colon {\Bbb N}^{\oplus 2} \to {\Bbb Z})$ with $\varphi(m, n) = m + n$ and $\alpha(m, n) = x^my^n$ forms a pre-log ring, with trivial actions.
\end{example}

As Example \ref{ex:comm} illustrates, there are strict inclusions \[{\rm CCPreLog} \subset {\rm CPreLog} \subset {\rm PreLog}\] involving classically commutative pre-log rings, pre-log rings $(A, \varphi \colon M \to G)$ with $A$, $M$, and $G$ commutative with trivial actions, and the category of all pre-log rings. 

\begin{example}[Logarithmic differential operators] Let us consider the triple \[({\Bbb Z}\langle t, \partial \rangle / (\partial t - t\partial - t), \gamma \colon {\Bbb N} \to {\Bbb Z}, \alpha),\] where $\alpha(n) = t^n$. Here, ${\Bbb Z}$ is assumed to act trivially on ${\Bbb N}$, while the action on ${\Bbb Z}\langle t, \partial \rangle / (\partial t - t\partial - t)$ is determined by $n.t = t$ and $n.\partial = \partial - n$. Observe that the defining relation $t\partial = (\partial - 1)t$ implies $t^n\partial = (\partial - n)t^n$, which is exactly (the interesting part of) the Peiffer condition. Thus, this defines a pre-log ring. 

This example is (an integral model of) the sheaf of logarithmic differential operators associated with ${\Bbb A}^1$ with log structure at $0$, see e.g.\ \cite[Example 3.3]{KT19}. The same principle applies to endow the sheaf of log differential operators on any log smooth log scheme with the structure of a sheaf of pre-log rings, using the explicit description of \cite[Corollary 3.2]{KT19}.  
\end{example}

\section{Associative log rings}\label{sec:asslog} Given a pre-log ring $(A, \varphi \colon M \to G)$, we say that it is \emph{log} if $\alpha^{-1}A^\times \to A^\times$ is an isomorphism. We shall write ${\rm Log} \subset {\rm PreLog}$ for the resulting full subcategory of pre-log rings that are log. In this section, we give an explicit construction of the left adjoint to this inclusion, which specializes to the usual logification construction (``passing to the associated log structure'') when applied to classically commutative pre-log rings. 

\subsection{Logification} Recall that, given a classically commutative pre-log ring $(A, M)$, its \emph{logification} is defined as $(A, M^a)$, where $M^a := M \oplus_{\alpha^{-1}A^\times} A^\times$, and its structure map $\alpha^a \colon M^a \to (A, \cdot)$ is induced by $\alpha$ and the inclusion $A^\times \subset (A, \cdot)$.

In particular, the pushout defining $M^a$ is here one of \emph{commutative} monoids. While it is tempting to define a logification construction of associative pre-log rings using a pushout of associative monoids, this is problematic for at least two reasons: First, this will seldom recover the ordinary logification construction when the input is classically commutative, and second, pushouts of associative monoids are generally unpleasant to work with. 

The following construction highlights how the additional data present in a general associative pre-log ring allow us to circumvent the lack of commutativity: 

\begin{construction}[Logification]\label{constr:logification} Given a pre-log ring $(A, \varphi \colon M \to G)$, we first consider the semidirect product $M \ltimes A^\times$ with $(m, u)(n, v) := (mn, u(\varphi(m).v))$, and similarly for $G \ltimes A^\times$. Then $\varphi$ extends to a map $\hat{\varphi} \colon M \ltimes A^\times \to G \ltimes A^\times$. As we explain in the proof of Lemma \ref{lem:logiprelim} below, this is a log crossed module with respect to the action \begin{equation}\label{eq:prelimcrossed}(g, u).(m, v) := (g.m, u(g.v)(\varphi(g.m).u^{-1}));\end{equation} this is the formula you obtain by combining the natural action $(g, 1).(m, v) := (g.m, g.v)$ with the conjugation-action of $(1, u)$. We obtain a pre-log ring \begin{equation}\label{eq:prelimlog}(A, \hat{\varphi} \colon M \ltimes A^\times \to G \ltimes A^\times, \hat{\alpha}),\end{equation} by Lemma \ref{lem:logiprelim} below, where $(g, u).a := {\rm Ad}_u(g.a) := u(g.a)u^{-1}$ and $\hat{\alpha}(m, u) := u\alpha(m)$. We define $M^a$ and $G^a$ to be the quotients \[M^a := (M \ltimes A^\times)/((1, \alpha(t)) \sim (t, 1)), \quad G^a := (G \ltimes A^\times)/((1, \alpha(t)) \sim (\varphi(t), 1))\] for $t \in \alpha^{-1}A^\times$. We now argue that the structures of \eqref{eq:prelimlog} descend to these quotients: It is clear that $\hat{\alpha}$ induces a map $\alpha^a \colon M^a \to (A, \cdot), \quad [m, u] \mapsto u\alpha(m)$. Next, we observe that the Peiffer condition implies that \[\varphi(t).a = \alpha(t)a\alpha(t)^{-1} = {\rm Ad}_{\alpha(t)}(a)\] for all $t \in \alpha^{-1}A^\times$. From this, we find that $(1, \alpha(t)).a = {\rm Ad}_{\alpha(t)}(a) = (\varphi(t), 1).a$ for all $t \in \alpha^{-1}A^\times$ and $a \in A$, so that this action gives rise to an action of $G^a$ on $A$ by ring automorphisms. By construction, $\hat{\varphi}$ gives a map $\varphi^a \colon M^a \to G^a$.

It remains to describe the action of $G^a$ on $M^a$. We take this to be induced from the corresponding action in \eqref{eq:prelimlog}: \[[g, u].[m, v] = [g.m, u(g.v)(\varphi(g.m).u^{-1})].\] In order to see that this is well-defined, let us first argue that this formula gives an action of $G \ltimes A^\times$ on $M^a$: Indeed, for any $t \in \alpha^{-1}A^\times$, we have \begin{equation}\label{eq:newaction1}(g, u).(t, 1) = (g.t, u(\varphi(g.t).u^{-1})) = (1, u)(g.t, 1)(1, u^{-1}),\end{equation} while \[(g, u).(1, \alpha(t)) = (1, u(g.\alpha(t))(\varphi(1).u^{-1})) = (1, u(g.\alpha(t))u^{-1}).\] Using equivariance of $\alpha$, we rewrite this as \begin{equation}\label{eq:newaction2}(1, u\alpha(g.t)u^{-1}) = (1, u)(1, \alpha(g.t))(1, u^{-1}).\end{equation} The expressions \eqref{eq:newaction1} and \eqref{eq:newaction2} are equivalent in $M^a$. To show that the resulting action of $G \ltimes A^\times$ on $M^a$ descends to an action of $G^a$, we first observe that \[[(1, \alpha(t)).(n, v)] = [(1, \alpha(t))(n, v)(1, \alpha(t)^{-1})] = [(t, 1)(n, v)(1, \alpha(t)^{-1})]\] for all $t \in \alpha^{-1}A^\times$. This is equal to $[(t, 1)(n, v)][(t, 1)]^{-1}$, which equals \[[((\varphi(t), 1).(n, v))(t, 1)][(t, 1)]^{-1} = [(\varphi(t), 1).(n, v)]\] by the Peiffer identity.
\end{construction}

\begin{example} Let $(A, M)$ be a classically commutative pre-log ring. Then the balancing relations defining $M^a$ are precisely those of the pushout $M \oplus_{\alpha^{-1}A^\times} A^\times$ of commutative monoids. 
\end{example}

\noindent The above construction enjoys the expected universal property:

\begin{proposition}\label{prop:logification} Let $(A, \varphi \colon M \to G)$ be a pre-log ring. The assignment \[(A, \varphi \colon M \to G) \mapsto (A, \varphi^a \colon M^a \to G^a)\] defines a functor ${\rm PreLog} \to {\rm Log}$ which is left adjoint to the inclusion. 
\end{proposition}

The proof of Proposition \ref{prop:logification} requires some preparation. Let us first verify that the auxiliary object of Construction \ref{constr:logification} is indeed a pre-log ring: 

\begin{lemma}\label{lem:logiprelim} The construction \eqref{eq:prelimlog} defines a pre-log ring. 
\end{lemma}

\begin{proof} Let us first show that the map \[\hat{\varphi} \colon M \ltimes A^\times \to G \ltimes A^\times, \quad (m, u) \mapsto (\varphi(m), u)\] defines a crossed semimodule with respect to the action \eqref{eq:prelimcrossed}. It is straightforward to check that it is an equivariant monoid map, while the Peiffer identity follows from that for the crossed semimodule $\varphi$: Indeed, $(\hat{\varphi}(m, u).(n, v))(m, u)$ equals \begin{align*} & ((\varphi(m).n, u(\varphi(m).v)(\varphi(\varphi(m).n).u^{-1})))(m, u) \\ =  & ((\varphi(m).n)m, u(\varphi(m).v)(\varphi(\varphi(m).n).u^{-1})(\varphi(\varphi(m).n).u)) \\ = &  (mn, u(\varphi(m).v)),\end{align*} where in the last equality we have used that $a \mapsto \varphi(\varphi(m).n).a$ is a ring map in the second coordinate and the Peiffer identity for $\varphi$ in the first. 

To see that $\hat{\alpha}$ is a monoid map, we use the Peiffer identity. Indeed, \[\hat{\alpha}((m, u)(n, v)) = u(\varphi(m).v)\alpha(m)\alpha(n) = u\alpha(m)v\alpha(n) = \hat{\alpha}(m, u)\hat{\alpha}(n, v).\] The map $\hat{\alpha}$ is readily seen to satisfy the conditions of a pre-log ring, using that $\alpha$ does so. 
\end{proof}

\begin{corollary} Let $(A, M)$ be a pre-log ring. Then the construction \[(A, \varphi^a \colon M^a \to G^a, \alpha^a)\] of Construction \ref{constr:logification} defines a pre-log ring. 
\end{corollary}

\begin{proof} By Lemma \ref{lem:logiprelim}, the intermediate object \eqref{eq:prelimlog} is a pre-log ring. By the factorizations verified in Construction \ref{constr:logification}, the relevant quotients thus participate in the asserted pre-log structure. 
\end{proof}

\begin{proof}[Proof of Proposition \ref{prop:logification}] The proposed logification is a log ring for the same reason that it is so in the classically commutative case: $\alpha^a_{|(\alpha^a)^{-1}A^\times} \colon (\alpha^a)^{-1}A^\times \to A^\times$ is an isomorphism with inverse $u \mapsto [1, u]$. There is a canonical map \[(i, i^\flat, i^\#) \colon (A, \varphi \colon M \to G) \to (A, \varphi^a \colon M^a \to G^a)\] of pre-log rings, factoring through the intermediate pre-log ring \eqref{eq:prelimlog}. Suppose that \[(f, f^\flat, f^\#) \colon (A, \varphi \colon M \to G) \to (B, \psi \colon N \to H)\] is a map of pre-log rings with $(B, \psi \colon N \to H)$ log. We shall argue that there is a unique extension $(\hat{f}, \hat{f}^\flat, \hat{f}^\#) \colon (A, \hat{\varphi} \colon M \ltimes A^\times \to G \ltimes A^\times) \to (B, \psi \colon N \to H)$ which admits a further unique extension through the logification. 

Since $(B, \psi \colon N \to H)$ is log, its structure map $\beta \colon N \to (B, \cdot)$ restricts to an isomorphism $N^\times \xrightarrow{\cong} B^\times$. For $v \in B^\times$, let $\beta^{-1}(v)$ denote the corresponding unit in $N$. Compatibility with the structure maps forces the formulas \[\hat{f}(a) = f(a), \quad \hat{f}^\flat(m, u) = \beta^{-1}(f(u))f^\flat(m), \quad \hat{f}^\#(g, u) = \psi(\beta^{-1}(f(u)))f^\#(g).\] This defines a morphism of pre-log rings, as is readily checked. 

It remains to see that this indeed extends uniquely over the logification: For $t \in \alpha^{-1}A^\times$, we have \[\hat{f}^\flat(1, \alpha(t)) = \beta^{-1}(f(\alpha(t))) = \beta^{-1}(\beta(f^\flat(t))) = f^\flat(t) = \hat{f}^\flat(t, 1),\] which also gives the necessary formula for descent to $G^a$ upon applying $\psi$.
\end{proof}

\subsection{Examples of log rings} These can now be obtained by logifying the examples of pre-log rings already discussed. For example, the adjoint pre-log ring \[(A, q \colon {\rm Ad}(A) \to {\rm Aut}(A))\] of Example \ref{ex:adaprelog} is already a log ring. 

\section{Logarithmic bimodules and the log enveloping algebra}\label{sec:logenv}

We now introduce the notion of a bimodule over a pre-log ring $(A, M)$, and we prove that they may be described as left modules over an associative ring $(A, M)^e$ which we refer to as the \emph{log enveloping algebra}. The technical core of this section is Proposition \ref{prop:logenvlogi}, in which we prove that this algebra is invariant under the logification construction of the previous section. 

\subsection{Logarithmic bimodules} Let $A$ be an associative ring. Recall that an $(A, A)$-bimodule consists of an abelian group $J$ that is simultaneously a left and a right $A$-module, subject to the compatibility condition $(aj)a' = a(ja')$. In the context of pre-log rings, we propose the following definition:

\begin{definition}\label{def:logbimod} Let $(A, \varphi \colon M \to G, \alpha)$ be a pre-log ring. An \emph{$(A, M)$-bimodule} consists of an $(A, A)$-bimodule $J$ that is simultaneously a left $G$-module, subject to the compatibility conditions \[g(ajb) = (g.a)gj(g.b) \quad \text{and} \quad \alpha(m)j\alpha(n) = \varphi(m)j\alpha(mn).\] We write ${\rm BiMod}(A, M)$ for the resulting category of log bimodules. 
\end{definition}

\subsection{The logarithmic enveloping algebra} The data of a bimodule over an associative ring $A$ is equivalent to a left module over the enveloping algebra $A^e$, defined as $A \otimes A^{\rm op}$. We now seek an analog of this observation in the context of pre-log rings. As briefly explained in the introduction, there is an obvious candidate in the context of classically commutative pre-log rings: namely the ring \[(A, M)^e := (A \otimes A) \otimes_{{\Bbb Z}[M \oplus M]} {\Bbb Z}[M^{\rm gp} \oplus M],\] where the relative tensor product is formed using the structure map of $(A, M)$ and \begin{equation}\label{eq:shearing}M \oplus M \xrightarrow{(m, m') \mapsto (\gamma(m), mm')} M^{\rm gp} \oplus M,\end{equation} where $\gamma \colon M \to M^{\rm gp}$ is the canonical map. There is by now very strong evidence that this is the correct analog of the enveloping algebra in the classical commutative setting: For example, by \cite[Section 13]{Rog09} or \cite[Proposition 1.4]{BLPO23a}, Rognes' log Hochschild homology ${\rm HH}_*(A, M)$ may be computed by the ${\rm Tor}$-groups ${\rm Tor}_*^{(A, M)^e}(A, A)$ under suitable flatness hypotheses. 

The map \eqref{eq:shearing} is a monoid map precisely when $M$ is commutative. Thus, for a general pre-log ring $(A, M)$, it is not immediately clear what should play the role of $(A, M)^e$. We now explain how the additional structure present in a pre-log ring allows us to construct an adequate replacement that recovers the above construction in the classical commutative case: 

\begin{construction}[Logarithmic enveloping algebra] Let $(A, M)$ be a pre-log ring, and let us consider the semidirect product $G \ltimes M^{\rm op}$ with multiplication given by $(g, m)(h, m') = (gh, (g.m')m)$. We now form the balanced tensor product (not the pushout) $(A, M)^e$ of the diagram \[A \otimes A^{\rm op} \xleftarrow{} {\Bbb Z}[M] \otimes {\Bbb Z}[M]^{\rm op} \xrightarrow{m \otimes m' \mapsto (\varphi(m), mm')} {\Bbb Z}[G \ltimes M^{\rm op}]\] of associative rings. In order to give this a ring structure, we first observe that $A \otimes A^{\rm op}$ admits an action of $G \ltimes M^{\rm op}$ by ring automorphisms, determined by \begin{equation}\label{eq:aaopaction}(g, m).(a \otimes b) := g.a \otimes (\varphi(m)^{-1}g).b.\end{equation} We equip $(A, M)^e$ with the multiplication \[((a \otimes b) \otimes (g, m))((a' \otimes b') \otimes (g', m')) := (a \otimes b)((g, m).(a' \otimes b')) \otimes (g, m)(g', m'),\] which, once spelled out, reads \[(a(g.a') \otimes ((\varphi(m)^{-1}g).b')b) \otimes (gg', (g.m')m).\] This gives $(A, M)^e$ the structure of a unital associative ring; for example, observe \begin{equation}\label{eq:logenvring1}((\alpha(m) \otimes \alpha(m')) \otimes (1, 1))((a \otimes b) \otimes (g, n)) = (\alpha(m)a \otimes b\alpha(m')) \otimes (g, n)\end{equation} By the Peiffer identities $\alpha(m)a = (\varphi(m).a)\alpha(m)$ and \[b\alpha(m') = (\varphi(m').(\varphi(m')^{-1}.b))\alpha(m') = \alpha(m')(\varphi(m')^{-1}.b),\] we may rewrite \eqref{eq:logenvring1} as \[((\varphi(m).a)\alpha(m) \otimes \alpha(m')(\varphi(m')^{-1}.b)) \otimes (g, n).\] Pulling out the factor $\alpha(m) \otimes \alpha(m')$ in the $A \otimes A^{\rm op}$-term and using the defining balancing relation, we find that this equals \[(\varphi(m).a \otimes \varphi(m')^{-1}.b) \otimes (\varphi(m), mm')(g, n) = ((1 \otimes 1) \otimes (\varphi(m), mm'))((a \otimes b) \otimes (g, n)).\] Thus, this twisted multiplication does respect the balancing relation of $(A, M)^e$. 
\end{construction}

\begin{example} If $(A, M)$ is a classically commutative pre-log ring, considered a pre-log ring with trivial actions, then $(A, M)^e$ is the pushout of the diagram \[A \otimes A \xleftarrow{} {\Bbb Z}[M \oplus M] \xrightarrow{(m, m') \mapsto (\gamma(m), mm')} {\Bbb Z}[M^{\rm gp} \oplus M]\] in the category of commutative rings. 
\end{example}

\begin{lemma}\label{lem:modbimod} There is a canonical equivalence of categories \[{\rm LMod}((A, M)^e) \simeq {\rm BiMod}(A, M)\] which equals the identity functor on underlying abelian groups.  
\end{lemma}

\begin{proof} Given a left $(A, M)^e$-module $J$, we endow $J$ with the $(A, M)$-bimodule structure determined by \[aj := ((a \otimes 1) \otimes (1, 1)).j, \quad ja := ((1 \otimes a) \otimes (1,1)).j, \quad gj := ((1 \otimes 1) \otimes (g, 1)).j.\] From this, the first condition of an $(A, M)$-bimodule is a direct consequence of the multiplication law in $(A, M)^e$, while its balancing relation gives the second condition. Conversely, we associate to an $(A, M)$-bimodule $J$ the left $(A, M)^e$-module determined by the formula $((a \otimes b) \otimes (g, m)).j := a(gj)\alpha(m)b$. Let us verify that this indeed defines an $(A, M)^e$-module structure: It respects the balancing relation, since $\alpha(m)j\alpha(n) = (\varphi(m)j)\alpha(m)\alpha(n)$ by the second defining relation of an $(A, M)$-bimodule. The associativity of the action is a straightforward (albeit slightly tedious) exercise in applying the pre-log ring axioms. 
\end{proof}

\subsection{Logification invariance of the log enveloping algebra} In the classically commutative setting, most invariants of pre-log rings, such as the module of differentials, the cotangent complex, and Hochschild homology, are invariant under logification. These are all consequences of the following more general fact:

\begin{proposition}\label{prop:logenvlogi} Let $(A, M)$ be a pre-log ring. The canonical map \[i \colon (A, M)^e \to (A, M^a)^e, \quad (a \otimes b) \otimes (g, m) \mapsto (a \otimes b) \otimes ([g, 1], [m, 1])\] is an isomorphism of associative rings. 
\end{proposition}

\begin{proof} Let us observe that for any $(A, M^a)$-bimodule $J$, we have equalities \begin{equation}\label{eq:conjlogeq}[g, u]j = [1, u]([g, 1]j) = u([g, 1]j)u^{-1},\end{equation} where in the last equality we have used the log bimodule axiom $\alpha(m)j\alpha(n) = \varphi(m)j\alpha(mn)$ with $m = [1, u], j = [g, 1]j$, and $n = [1, u]^{-1}$. It follows that restriction of scalars ${\rm BiMod}(A, M^a) \to {\rm BiMod}(A, M)$ is an equivalence of categories. Indeed, for any $(A, M)$-bimodule we obtain an $(A, M^a)$-bimodule by declaring that $[g, u]j = u(gj)u^{-1}$; one checks that this satisfies the log bimodule conditions using the defining properties of pre-log rings. It is a well-defined inverse by \eqref{eq:conjlogeq}. By Lemma \ref{lem:modbimod}, restriction of scalars along $(A, M) \to (A, M^a)$ is an equivalence ${\rm LMod}((A, M^a)^e) \to {\rm LMod}((A, M)^e)$, thus is also extension of scalars. It follows that the adjunction unit $J \to (A, M^a)^e \otimes_{(A, M)^e} J$ is an isomorphism for any $(A, M)^e$-module $J$. The desired result follows from setting $J = (A, M)^e$.
\end{proof}

The following remark is probably only of interest to experts on log (topological) Hochschild homology: 

\begin{remark}\label{rem:logenvlogi} The proofs of logification invariance of (topological) Hochschild homology currently available in the literature follow one of the following strategies:
    
    \begin{enumerate}
        \item In the context of classically commutative pre-log rings, one can use the logarithmic ${\rm HKR}$-filtration of \cite[Theorem 1.1]{BLPO23a} to reduce logification invariance of log Hochschild homology ${\rm HH}(A, M)$ to the log cotangent complex ${\Bbb L}_{(A, M)}$, which is the content of \cite[Theorem 8.16]{Ols05}.
        \item In the context of the log ring spectra developed in \cites{Rog09, RSS15, RSS18, RSS26}, there is instead a more technically involved argument using the Bousfield--Friedlander theorem \cite[Theorem B.4]{BF78}, see \cite[Theorem 4.24]{RSS15}. At a certain point in the argument, one uses that classifying spaces of groups are path-connected. 
    \end{enumerate}
\noindent While the first strategy above reduces from log Hochschild homology to the cotangent complex, it is also possible to argue in the opposite direction. For simplicity, let us consider a classically commutative pre-log ring $(A, M)$ for which $A$ is ${\Bbb Z}[M]$-flat (or simply derive everything in sight), and consider the \emph{pointed} simplicial tensors $S^n \odot_A (A, M)^e$ in simplicial augmented $A$-algebras. For $n = 1$, this is ${\rm HH}(A, M)$ by \cite[Proposition 5.8]{Lun21}, while as $n$ tends to $\infty$ (stabilization), it is ${\Bbb L}_{(A, M)}$ by \cite[Section 9]{Lun21} or \cite[Section 13]{Rog09}. From this perspective, the path-connectivity of classifying spaces used in (2) translates to path-connectivity of $S^1$. For this reason, the authors did not expect $S^0 \odot_A (A, M)^e \cong (A, M)^e$ to be logification invariant in general. But this is precisely the content of Proposition \ref{prop:logenvlogi}, and logification invariance of the logarithmic differentials, cotangent complex, and Hochschild homology now become formal consequences, at least when $A$ is ${\Bbb Z}[M]$-flat. It also seems plausible that an analogous statement is true in the spectral context of e.g.\ \cite{RSS26}, but we did not pursue this. 
\end{remark}

\section{Properties of pre-log rings and their morphisms}\label{sec:properties} We now pause to record various analogs of definitions from classical log geometry, and we explain how they relate to each other. This is all heavily inspired by classical properties as for instance recorded in Ogus' book \cite{Ogu18}: The general principle is that the role of the group completion map $\gamma \colon M \to M^{\rm gp}$ in the classically commutative case is now played by the structure map $\varphi \colon M \to G$. 

\subsection{Integral pre-log rings} We say that a pre-log ring $(A, \varphi \colon M \to G)$ is \emph{integral} if $\varphi$ is injective. In the classically commutative case, this translates to $\gamma \colon M \to M^{\rm gp}$ being injective, that is, $M$ being cancellative. 

\subsection{Exact morphisms of pre-log rings} We say that a morphism \[(p, p^\flat, p^\#) \colon (B, \psi \colon N \to H) \to (A, \varphi \colon M \to G)\] is \emph{exact} if the canonical map $N \to M \times_G H$ is an isomorphism. 

\begin{remark} In analogy with the terminology used in \cite{Rog09}, we are tempted to say that an exact morphism $(p, p^\flat, p^\#)$ is \emph{replete} if $p^\#$ is moreover surjective. While we refrain from using this terminology, this condition will certainly play a role in our forthcoming results (cf.\ Lemma \ref{lem:exstrict}). 
\end{remark}
 
\subsection{Strict morphisms of log rings}\label{subsec:logsqzero} A map $(B, \psi \colon N \to H) \to (A, \varphi \colon M \to G)$ of log rings is \emph{strict} if it induces an isomorphism \[(A, \psi \colon N \to H)^a \xrightarrow{\cong} (A, \varphi \colon M \to G),\] where the source is the logification of the pre-log structure on $A$ with structure map $N \to M \to (A, \cdot)$ and $H$-action restricted along $H \to G$.   

\subsection{Log square-zero extensions} We follow the convention that square-zero extensions of rings are surjective. With strictness in place, we define:

\begin{definition} A morphism $(p, p^\flat, p^\#) \colon (B, \psi \colon N \to H) \to (A, \varphi \colon M \to G)$ of log rings is \emph{square-zero} if it is strict and $p$ is square-zero. 
\end{definition}

\begin{lemma}\label{lem:psharpsurj} Let $(p, p^\flat, p^\#) \colon (B, \psi \colon N \to H) \to (A, \varphi \colon M \to G)$ be a square-zero extension of log rings. Then $p^\flat$ and $p^\#$ are surjective. 
\end{lemma}

\begin{proof} We first argue that $p^\#$ is surjective. Since $(p, p^\flat, p^\#)$ is strict, $G$ is isomorphic to the underlying group of the logification of $(A, \psi \colon N \to H)$. Let us write $[h, u]$ for an element of $G$. Now we use that $p$ being square-zero implies that the restriction $p_{|B^\times} \colon B^\times \to A^\times$ is surjective, as well as $(B, \psi \colon N \to H)$ being log implying that the structure map $\beta \colon N \to (B, \cdot)$ restricts to an isomorphism $\beta_{|N^\times} \colon N^\times \xrightarrow{\cong} B^\times$, to infer that there exists $n \in N^\times$ such that $p(\beta(n)) = u$. The balancing relation in $G$, that is, the group part of the logification of $(A, \psi \colon N \to H)$, now gives \[[h, u] = [1, u][h, 1] = [\psi(n), 1][h, 1] = [\psi(n)h, 1],\] which equals $p^{\#}(\psi(n)h)$.  Similarly, by the strictness hypothesis, we can write any $m \in M$ as $[n, u]$ for $u \in A^\times$ and $n \in N$. Choosing $v \in B^\times$ such that $p(v) = u$, we have $[n, u] = [n, p(v)] = [\beta^{-1}(v)n, 1]$, so that $m = p^\flat(\beta^{-1}(v)n)$.
\end{proof}

In ordinary algebra, the property of $p \colon B \to A$ being square-zero equips its kernel $J$ with an $(A, A)$-bimodule structure in a canonical way. As further evidence that our proposed notion of log bimodules is sound, we have:

\begin{lemma}\label{lem:sqzerobimod1} Let $(p, p^\flat, p^\#) \colon (B, \psi \colon N \to H) \to (A, \varphi \colon M \to G)$ be a square-zero extension and let $J := {\rm ker}(p)$. Then the actions \[p(b)j := bj, \quad jp(b) := jb, \quad p^\#(h)j := h.j\] are well-defined and exhibit $J$ as an $(A, M)$-bimodule.   
\end{lemma}

\begin{proof} Since both $p$ and $p^\#$ are surjective by Lemma \ref{lem:psharpsurj}, we have indeed specified a proposed action from all of $A$ and $G$. It is classical that this endows $J$ with an $(A, A)$-bimodule structure. Since $p$ is equivariant, $h.j$ is an element of $J$. 

We will first argue that $J$ is an $(A, N)$-bimodule. The first condition $h.(bjb') = (h.b)(h.j)(h.b')$ simply follows from $H$ acting by ring (auto)morphisms. The second condition reduces to the Peiffer identity: Indeed, we have that \[p(\beta(n))jp(\beta(n')) = \beta(n)j\beta(n') = (\psi(n).j)\beta(nn') = (\psi(n).j)p(\beta(nn')).\] By Lemma \ref{lem:modbimod}, we have specified $J$ as a left $(A, N)^e$-module. By strictness and Proposition \ref{prop:logenvlogi}, the canonical map $(A, N)^e \to (A, M)^e$ is an isomorphism, so this gives $J$ the desired $(A, M)^e$-module structure. 
\end{proof}

In light of Lemma \ref{lem:sqzerobimod1}, we may write e.g.\ ``Let $(p, p^\flat, p^\#) \colon (B, \psi \colon N \to H) \to (A, \varphi \colon M \to G)$ be a square-zero extension by an $(A, M)^e$-module $J$.'' In these cases, it is always the $(A, M)^e$-module structure stemming from Lemma \ref{lem:sqzerobimod1} that is meant. 

\subsection{The structure of square-zero extensions} We now proceed to study some associative analogs of the results in \cite[Section 4]{SSV16}, which will be used in the proof of Theorem \ref{thm:pullback}. We have the following analog of \cite[Lemma 4.15]{SSV16}. 

\begin{proposition}\label{prop:strictsqzeroexact} Let $(p, p^\flat, p^\#) \colon (B, \psi \colon N \to H) \to (A, \varphi \colon M \to G)$ be a strict square-zero extension of integral log rings. Then $(p, p^\flat, p^\#)$ is exact. 
\end{proposition}

\noindent The proof of Proposition \ref{prop:strictsqzeroexact} requires some preparation:

\begin{lemma}\label{lem:strictsquarezeroexact2} With notation as in Proposition \ref{prop:strictsqzeroexact}, there is an equality ${\rm ker}(p^\#) = \psi({\rm ker}(p^\flat))$ of subgroups of $H$.  
\end{lemma}

\begin{proof} The square-zero and log conditions imply that ${\rm ker}(p^\flat)$ is a subgroup of $N^\times$. By equivariance of $p^\flat$, it is stable under the action of $H$, and by equivariance of $\psi$, the image $\psi({\rm ker}(p^\flat))$ is a normal subgroup of $H$. Spelling out the strictness condition in this case shows that $p^\#$ induces an isomorphism $H / \psi({\rm ker}(p^\flat)) \cong G$, which gives the result. 
\end{proof}

\begin{proof}[Proof of Proposition \ref{prop:strictsqzeroexact}] The claim is that the canonical map \begin{equation}\label{eq:strictcanmap}N \to M \times_G H, n \mapsto (p^\flat(n), \psi(n))\end{equation} is an isomorphism. It is injective since $\psi$ is, and so it suffices to prove that it is surjective. Given $(m, h) \in M \times_G H$, we can by Lemma \ref{lem:psharpsurj} choose $n \in N$ with $p^\flat(n) = m$. Since $(p^\flat(n), h) \in M \times_G H$, we have that $h\psi(n)^{-1} \in {\rm ker}(p^\#)$, which equals $\psi({\rm ker}(p^\flat))$ by Lemma \ref{lem:strictsquarezeroexact2}. So choose $t \in {\rm ker}(p^\flat)$ so that $\psi(t) = h\psi(n)^{-1}$. Then $tn$ hits $(m, h)$ under \eqref{eq:strictcanmap}, which concludes the proof. 
\end{proof}

\noindent Conversely, we have the following analog of \cite[Proposition I.4.2.1(5)]{Ogu18}:

\begin{lemma}\label{lem:exstrict} Let $(p, p^\flat, p^\#) \colon (B, \psi \colon N \to H) \to (A, \varphi \colon M \to G)$ be an exact morphism of log rings. If $p^\#$ is surjective, then $(p, p^\flat, p^\#)$ is strict. 
\end{lemma}

\begin{proof} We show that $(A, \varphi \colon M \to G)$ enjoys the universal property of the logification of $(A, \psi \colon N \to H)$. So let $(f, f^\flat, f^\#) \colon (A, \psi \colon N \to H) \to (C, \rho \colon K \to L)$ be a morphism with $(C, \rho \colon K \to L)$ log. By exactness, we have $N \cong M \times_G H$. If $h \in {\rm ker}(p^\#)$, then $(1_M, h) \in N$ satisfies $\alpha_C(f^\flat(1_M, h)) = f(\alpha_A(p^\flat(1_M, h))) = 1_C$. Since $(C, \rho \colon K \to L)$ is log, this implies that $f^\flat(1_M, h) = 1_K$ and thus that $f^\#(h) = 1_L$. Thus, given elements $(m, h)$ and $(m, h') \in N$, we have $f^\flat(m, h) = f^\flat(m, h')$. 

In total, the map \[(f, \overline{f}^\flat, \overline{f}^\#) \colon (A, \varphi \colon M \to G) \to (C, \rho \colon K \to L)\] determined by $\overline{f}^\flat(m) = f^\flat(m, h_{\varphi(m)})$ and $ \overline{f}^\#(g) = f^\#(h_g)$, where e.g.\ $h_g$ is any lift of $g$ under $p^\#$, is a well-defined and necessarily unique factorization of $(f, f^\flat, f^\#)$. Since $(A, \varphi \colon M \to G)$ is log, this concludes the proof. 
\end{proof}

\section{Logarithmic bimodules as abelian group objects}\label{sec:abgrp}

We now work towards the proof of Theorem \ref{thm:abgrpob} from the introduction. As a first step, we establish the asserted equivalence of categories in Corollary \ref{cor:abgrpob}.

In the proof of the classically commutative analog of Theorem \ref{thm:abgrpob}, the key observation is that for fixed $(A, M)$, the category of log rings $(A, M) \to (B, N) \to (A, M)$ over and under $(A, M)$ with $(B, N) \to (A, M)$ exact no longer depends upon the log structure: The forgetful functor ${\rm CCLog}^{\rm ex}_{(A, M) // (A, M)} \to {\rm CRing}_{A // A}$ is an equivalence of categories (cf.\ \cite[Corollary 2.14]{Lun25}). As we shall see, the analogous fact is \emph{not} true in the general associative setting:

\subsection{Split exact log rings}\label{subsec:splitexact} To begin, we first prove an analog of \cite[Lemma 2.12]{Lun25}. As in \emph{loc.\ cit.}, consider an augmented pre-log ring \begin{equation}\label{eq:retractive}(A, \varphi \colon M \to G) \xrightarrow{(i, i^\flat, i^\#)} (B, \psi \colon N \to H) \xrightarrow{(p, p^\flat, p^\#)} (A, \varphi \colon M \to G).\end{equation} Let us write ${\rm ker}(p^\times) := {\rm ker}(B^\times \xrightarrow{p_{|B^\times}} A^\times)$. The splitting $A^\times \ltimes {\rm ker}(p^\times) \xrightarrow{\cong} B^\times$, where $A^\times$ acts on ${\rm ker}(p^\times)$ by conjugation, induces the pre-log structure \begin{equation}\label{eq:lognormalform}(B, \varphi_B \colon M \ltimes {\rm ker}(p^\times) \to G \ltimes {\rm ker}(p^\times)),\end{equation} where the formulas \[(m, u)(m', u') := (mm', u(i^\#(\varphi(m)).u')), \quad (g, u)(g', u') := (gg', u(i^\#(g).u'))\] give the multiplications on the source and target in $\varphi_B \colon M \ltimes {\rm ker}(p^\times) \to G \ltimes {\rm ker}(p^\times)$, which is a log crossed module with respect to the action \[(g, u).(m, v) := (g.m, (ui^\#(g).v)(i^\#(\varphi(g.m)).u^{-1})).\] More conceptually, $G$ acts on ${\rm ker}(p^\times)$ via automorphisms, hence so does $M$ via $\varphi$, and the resulting action participating in the log crossed module structure combines the natural actions \[(g, 1).(m, v) = (g.m, i^\#(g).v), \quad (1, u).(m, v) := (1, u)(m, v)(1, u^{-1}).\] This is conveniently similar to the log crossed module action involved in the logification construction, cf.\ \eqref{eq:prelimcrossed}:

\begin{lemma}\label{lem:retractivelog} In \eqref{eq:retractive}, assume that $(p, p^\flat, p^\#)$ is exact. Then $(B, \psi^a \colon N^a \to H^a)$ is isomorphic to the log ring \begin{equation}\label{eq:retractivelog}(B, \varphi_B^a \colon M^a \ltimes {\rm ker}(p^\times) \to G^a \ltimes {\rm ker}(p^\times))\end{equation} over and under $(A, \varphi^a \colon M^a \to G^a)$, with structure maps as in \eqref{eq:lognormalform}. 
\end{lemma}

\begin{proof} The exactness hypothesis implies that ${\rm ker}(p^\#) \cong {\rm ker}(p^\flat) := (p^{\flat})^{-1}(1)$, while the fact that $(p, p^\flat, p^\#)$ is split implies that there is an isomorphism \[(B, \varphi \ltimes {\rm id} \colon M \ltimes {\rm ker}(p^\flat) \to G \ltimes {\rm ker}(p^\flat)) \xrightarrow{\cong} (B, \psi \colon N \to H)\] of pre-log rings over and under $(A, \varphi \colon M \to G)$. This makes it transparent to prove that \eqref{eq:retractivelog} enjoys the universal property of the logification of $(B, \psi \colon N \to H)$: We must show that, given any diagram \[\begin{tikzcd}[column sep = tiny](A, \varphi \colon M \to G) \ar{d} \ar{r} & (A, \varphi^a \colon M^a \to G^a) \ar{d} \\ (B, \varphi \ltimes {\rm id} \colon M \ltimes {\rm ker}(p^\flat) \to G \ltimes {\rm ker}(p^\flat)) \ar{r} \ar{d} & (C,\xi \colon K \to L) \ar{d} \\ (A, \varphi \colon M \to G) \ar{r} & (A, \varphi^a \colon M^a \to G^a) \end{tikzcd} \] with $(C, \xi \colon K \to L)$ log, the middle horizontal arrow factors compatibly through \eqref{eq:retractivelog}. The extensions to $M^a$ and $G^a$ are determined by the structure map from the log ring $(A, \varphi^a \colon M^a \to G^a)$, while the extension to ${\rm ker}(p^\times)$ comes from the observation that ${\rm ker}(p^\flat) \cong {\rm ker}(p^\#)$ is a group: Thus, the composition ${\rm ker}(p^\flat) \to N \to (B, \cdot)$ factors through $B^\times$. Hence the extension from ${\rm ker}(p^\flat)$ to ${\rm ker}(p^\times)$ is forced to be the composition ${\rm ker}(p^\times) \to B^\times \to C^\times \cong K^\times \subset K \xrightarrow{\xi} L$, where the isomorphism uses that $(C, \xi \colon K \to L)$ is log.  
\end{proof}

For a log crossed module $\varphi \colon M \to G$, let us write ${\rm PreLog}^{\varphi}$ for the category of pre-log rings whose log crossed module is $\varphi$. The following is the associative analog of \cite[Corollary 2.14]{Lun25}:

\begin{lemma}\label{lem:forgetfulphi} Let $(A, \varphi \colon M \to G)$ be a log ring. The forgetful functor \[{\rm Log}^{\rm ex}_{(A, M) // (A, M)} \xrightarrow{} {\rm PreLog}^{\varphi}_{(A, M) // (A, M)}, \quad (B, \psi \colon N \to H) \mapsto (B, \varphi \colon M \to G)\] is an equivalence of categories. 
\end{lemma}

\begin{proof} We claim that the inverse is given by logification. This boils down to the claim that for any $(B, \psi \colon N \to H) \in {\rm Log}^{\rm ex}_{(A, M) // (A, M)}$, the induced map \[(B, \varphi \colon M \to G)^a \to (B, \psi \colon N \to H)\] is an isomorphism. By Lemma \ref{lem:retractivelog}, both are compatibly isomorphic to \eqref{eq:retractivelog}, which concludes the proof. 
\end{proof}

\begin{lemma}\label{lem:abphi} Let $(A, \varphi \colon M \to G)$ be a pre-log ring. There is an equivalence \[{\rm Ab}({\rm PreLog}^\varphi_{(A, M) // (A, M)}) \simeq {\rm BiMod}{(A, M)}\] relating abelian group objects in ${\rm PreLog}^\varphi_{(A, M) // (A, M)}$ with $(A, M)$-bimodules. 
\end{lemma}

\begin{proof} The forgetful functor ${\rm PreLog}_{(A, M)//(A, M)}^\varphi \to {\rm Ring}_{A//A}$ preserves products, so it induces a functor on abelian group objects. Hence, the ring component of an abelian group object $(B, \varphi \colon M \to G)$ in ${\rm PreLog}_{(A, M)//(A, M)}^\varphi$ is of the form $A \oplus J$ for an $(A, A)$-bimodule $J$. The proposed $G$-action is defined via the formula \[g.(a, j) = (g.a, gj),\] which we now proceed to prove determines an $(A, M)$-bimodule structure on $J$. The first log bimodule axiom $g(ajb) = (g.a)gj(g.b)$ follows from \[g.(0, ajb) = g.((a, 0)(0, j)(b, 0)) = (g.a, 0)(0, gj)(g.b, 0).\] The second log bimodule axiom $\alpha(m)j\alpha(n) = (\varphi(m)j)\alpha(mn)$ follows from \[(\alpha(m), 0)(0, j) = \varphi(m).(0, j)(\alpha(m), 0) = (0, \varphi(m)j)(\alpha(m), 0) = (0, \varphi(m)j\alpha(m)),\] where we have used the Peiffer identity of $(B, \varphi \colon M \to G)$ in the first equality, and multiplied by $\alpha(n)$ from the right in the resulting equality $\alpha(m)j = \varphi(m)j\alpha(m)$.  

Conversely, an $(A, M)$-bimodule $J$ determines an abelian group object \begin{equation}\label{eq:abgrpob}(A \oplus J, \varphi \colon M \to G)\end{equation} in ${\rm PreLog}^{\varphi}_{(A, M) // (A, M)}$ by setting $g.(a, j) := (g.a, gj)$. Combining the pre-log axiom for $(A, M)$ and the log bimodule axiom for $J$, we obtain the equality \[(\alpha(m), 0)(a, j) = (\alpha(m)a, \alpha(m)j) = ((\varphi(m).a)\alpha(m), \varphi(m)j\alpha(m)),\] which equals $\varphi(m).(a, j)(\alpha(m), 0)$, which gives the desired pre-log axiom for the abelian group object $(A \oplus J, \varphi \colon M \to G)$. 
\end{proof}

\noindent This suffices to establish the equivalence of categories of Theorem \ref{thm:abgrpob}:

\begin{corollary}\label{cor:abgrpob} Let $(A, \varphi \colon M \to G)$ be a log ring. There is an equivalence of categories ${\rm Ab}({\rm Log}^{\rm ex}_{(A, M) // (A, M)}) \simeq {\rm BiMod}(A, M)$. 
\end{corollary}

\begin{proof} This follows from combining Lemmas \ref{lem:forgetfulphi} and \ref{lem:abphi}. 
\end{proof}

\begin{remark}\label{rem:diffalt} One natural variation of the category ${\rm PreLog}$ is the full subcategory ${\rm PreLog}^{\rm gp}$ consisting of those pre-log rings $(A, \varphi \colon M \to G)$ with the property that the induced map $M^{\rm gp} \to G$ is an isomorphism. For the purposes of developing the basic theory here, and in particular proving Corollary \ref{cor:abgrpob} in its stated form, this is a perfectly fine option. Given that allowing for non-trivial actions and the Peiffer identities remain essential also in this case, we find it more natural to build the theory for general $G$. 
\end{remark}

\section{Logarithmic derivations and Hochschild cohomology}\label{sec:loghch} We now proceed to study the natural definition of logarithmic Hochschild cohomology ${\rm HH}^*((A, M) ; J)$ that arises from our construction of the log enveloping algebra $(A, M)^e$. After identifying this group in degree $0$, we generalize the classical notion of a log derivation to the context of associative pre-log rings and prove that it participates in the expected isomorphism with ${\rm HH}^1((A, M) ; J)$. The theory developed along the way will allow us to conclude the proof of Theorem \ref{thm:abgrpob} towards the end of the section. 

\subsection{Logarithmic Hochschild cohomology} Let $(A, \varphi \colon M \to G)$ be a pre-log ring and let $J$ be a left $(A, M)^e$-module. We define the \emph{log Hochschild cohomology} ${\rm HH}^*((A, M) ; J)$ of $(A, M)$ with coefficients in $J$ to be the ${\rm Ext}$-groups ${\rm Ext}^*_{(A, M)^e}(A, J)$. As an immediate consequence of Proposition \ref{prop:logenvlogi}, we obtain the following invariance: 

\begin{corollary} Let $(A, \varphi \colon M \to G)$ be a pre-log ring and let $J$ be an $(A, M)$-bimodule. There is an isomorphism ${\rm HH}^*((A, M^a) ; J) \xrightarrow{\cong} {\rm HH}^*((A, M) ; J)$. \qed
\end{corollary}

\begin{remark} In the first-named author's master's thesis \cite{vdH26}, the cochain complex obtained by gluing the standard complexes computing Hochschild and group cohomology along the log structure is considered and pursued as a candidate definition of log Hochschild \emph{co}homology. This produces the same groups as ${\rm Ext}^*_{(A, M)^e}(A, J)$ for $* = 0, 1$. In degree $2$, one discovers a notion of square-zero extensions involving log rings that are potentially non-commutative in the most naive sense (cf.\ \cite[Definition 5.10, Remark 5.12]{Rog09}). With this definition, we would not know how to carry out most of the constructions pursued in this paper. 
\end{remark}

\subsection{The $G$-invariant center of a log bimodule} For $J$ a left $(A, M)^e$-module, we shall write $Z_A(J)^G := \{j \in J \mid aj = ja \text{ and } gj = j \text{ for all } a \in A, g \in G\}$ for the $G$-invariants of its centre. The following is analogous to ${\rm HH}^0(A ; J) \cong Z_A(J)$:

\begin{lemma}\label{lem:loghhzero} There is an isomorphism \[{\rm HH}^0((A, M); J) \xrightarrow{\cong} Z_A(J)^G,\] natural in the pre-log ring $(A, M)$ and $(A, M)^e$-module $J$. 
\end{lemma}

\begin{proof} Evaluating $f \in {\rm Ext}^0_{(A, M)^e}(A, J) \cong {\rm Hom}_{(A, M)^e}(A, J)$ at $1$ gives the isomorphism; the $A \otimes A^{\rm op}$-linearity ensures that this lands in $Z_A(J)$, while the ${\Bbb Z}[G]$-linearity ensures that we land in $Z_A(J)^G$. The inverse is determined by $j \mapsto (1 \mapsto j)$; this gives an $(A, M)^e$-linear map precisely because $j \in Z_A(J)^G$.
\end{proof}

\subsection{Logarithmic derivations} In degree $1$, log Hochschild cohomology captures a suitable notion of log derivations: 

\begin{definition}\label{def:logder} Let $(A, M)$ be a pre-log ring. A \emph{log derivation} $(d, \delta) \colon (A, M) \to J$ with values in a left $(A, M)^e$-module $J$ consists of a derivation $d \colon A \to J$ (i.e.\ $d(ab) = ad(b) + d(a)b$) and a crossed homomorphism $\delta \colon G \to J$ (i.e.\  $\delta(gh) = \delta(g) + g\delta(h)$) satisfying the following compatibility conditions:

\begin{enumerate}
    \item $d(g.a) = gd(a) + \delta(g)(g.a) - (g.a)\delta(g)$; and
    \item $d(\alpha(m)) = \delta(\varphi(m))\alpha(m)$
\end{enumerate}

\noindent for all $a \in A, g \in G$, and $m \in M$. We write ${\rm LogDer}((A, M), J)$ for the log derivations of $(A, M)$ with values in $J$.
\end{definition}

If $(A, M)$ is a classically commutative pre-log ring and $J$ is a symmetric bimodule with trivial action, then the first condition simply reads $d(a) = d(a)$, while the second reads $d(\alpha(m)) = \alpha(m)\delta(m)$. Thus, the above notion recovers the classical notion of logarithmic derivations in this case. In the general situation, fixed elements $j \in J$ always give rise to examples $(d_j, \delta_j)$ by declaring that $d_j(a) = ja - aj$ and $\delta_j(g) = j - gj$; we call these \emph{inner log derivations}. We write ${\rm InnDer}((A, M), J)$ for the inner log derivations of $(A, M)$ with values in $J$. 

Let us write $I_{(A, M)}$ for the kernel of $(A, M)^e \to A, (a \otimes b) \otimes (g, m) \mapsto a \alpha(m) b$. There is a logarithmic derivation with values in $I_{(A, M)}$ determined by the formulas \[D(a) := (1 \otimes a) \otimes (1, 1) - (a \otimes 1) \otimes (1, 1), \partial(g) := (1 \otimes 1) \otimes (1, 1) - (1 \otimes 1) \otimes (g, 1).\] Indeed, $D$ is obtained by composing the universal derivation with the canonical map from ${\rm ker}(A \otimes A^{\rm op} \to A)$, while $\partial$ is obtained by composing the universal crossed homomorphism with the canonical map from ${\rm ker}({\Bbb Z}[G] \to {\Bbb Z})$. The first condition of Definition \ref{def:logder} is the observation that, in the expression \[D(g.a) = (1 \otimes g.a) \otimes (1, 1) - (g.a \otimes 1) \otimes (1, 1)\] one may subtract $(1 \otimes g.a) \otimes (g, 1)$ from the first term to obtain $\partial(g)(g.a)$ and subtract $(g.a \otimes 1) \otimes (g, 1)$ from the second to obtain $(g.a)\partial(g)$. The difference between the subtracted terms is $gD(a)$, so the first condition of Definition \ref{def:logder} holds. The second condition is a consequence of the balancing relation \[(\alpha(m) \otimes 1) \otimes (1, 1) = (1 \otimes 1) \otimes (\varphi(m), m) = (1 \otimes \alpha(m)) \otimes (\varphi(m), 1)\] in the log enveloping algebra $(A, M)^e$. This is the universal log derivation:

\begin{lemma}\label{lem:unilogder} Let $(A, M)$ be a pre-log ring and $J$ an $(A, M)$-bimodule. The map \begin{equation}\label{eq:unilogder}{\rm Hom}_{{\rm LMod}((A, M)^e)}(I_{(A, M)}, J) \xrightarrow{f \mapsto (f \circ D, f \circ \partial)} {\rm LogDer}((A, M), J)\end{equation} is an isomorphism, natural in $(A, M)$ and $J$. 
\end{lemma}

\begin{proof} Beginning from the balancing relation \[(a \otimes b) \otimes (g, m) = (a \otimes \alpha(m)b) \otimes (g, 1),\] one deduces from the identities \[a \otimes \alpha(m)b = 1 \otimes a\alpha(m)b - (1 \otimes a - a \otimes 1)\alpha(m)b\] and $g = 1 - (1 - g)$ that \[(a \otimes b) \otimes (g, m) = (1 \otimes a \alpha(m) b) \otimes (1, 1) - D(a)\alpha(m)b - a\partial(g)\alpha(m)b.\] In particular, if $(a \otimes b) \otimes (g, m) \in I_{(A, M)}$, then this equals $-D(a)\alpha(m)b - a\partial(g)\alpha(m)b$, so that $I_{(A, M)}$ is generated by $D(a)$ and $\partial(g)$ for $a \in A, g \in G$. Thus, \eqref{eq:unilogder} is injective. Conversely, given a log derivation $(d, \delta)$, we have the map \[I_{(A, M)} \to J, \quad (a \otimes b) \otimes (g, m) \mapsto -d(a)\alpha(m)b - a\delta(g)\alpha(m)b,\] which provides a well-defined inverse of \eqref{eq:unilogder}. 
\end{proof}

In the non-logarithmic setting, derivations are closely related to (and sometimes defined in terms of) split square-zero extensions by means of the natural isomorphism ${\rm Der}(A, J) \cong {\rm Hom}_{{\rm Ring}_{/A}}(A, A \oplus J)$. In the context of classically commutative pre-log rings, there is an analogous relationship: With the language used here, the most convenient reference is \cite[Section 4]{Rog09}. In order to extend this to the context of associative pre-log rings, we first construct split square-zero extensions:

\begin{construction}[Split square-zero extensions]\label{constr:splitsqzero} Let $(A, \varphi \colon M \to G)$ be a pre-log ring and let $J$ be an $(A, M)$-bimodule. Consider the pre-log ring \[(A \oplus J, \varphi \ltimes J \colon M \ltimes J \to G \ltimes J, \alpha \ltimes J)\] with the following structure maps: $A \oplus J$ is the usual square-zero extension, $G \ltimes J$ is the usual semidirect product, and $M \ltimes J := M \times_G (G \ltimes J)$. We define the log crossed module by $(\varphi \ltimes J)(m, j) := (\varphi(m), j)$, which satisfies the axioms with the action $(g, j).(m, j') := (g.m, j + gj' - \varphi(g.m)j)$. The $G \ltimes J$-action on $A \oplus J$ is determined by the formula \[(g, j).(a, j') := {\rm Ad}_{(1, j)}(g.a, gj') := (g.a, gj' + j(g.a) - (g.a)j).\] The structure map $\alpha \ltimes J$ sends $(m, j)$ to $(\alpha(m), j\alpha(m))$.
\end{construction}

\noindent This construction enjoys the desired relationship with log derivations:

\begin{lemma}\label{lem:logdersqzero} Let $(A, M)$ be a pre-log ring and let $J$ be an $(A, M)^e$-module. There is an isomorphism \[{\rm LogDer}((A, M), J) \cong {\rm Hom}_{{\rm PreLog}_{/(A, M)}}((A, M), (A \oplus J, M \ltimes J))\] natural in $(A, M)$ and $J$.
\end{lemma}

\begin{proof} The right-hand ${\rm Hom}$-set sits as the pullback of \[\begin{tikzcd}[row sep = small]{\rm Hom}_{{\rm LogCrMod}_{/(M \to G)}}(M \to G, M \ltimes J \to G \ltimes J) \ar{d} \\  {\rm Hom}_{{\rm Grp}_{/G}}(G, G \ltimes J) \times {\rm Hom}_{{\rm Mon}_{/(A, \cdot)}}(M, (A \oplus J, \cdot)) \\ {\rm Hom}_{{\rm RingRep}_{/(A, G)}}((A, G), (A \oplus J, G \ltimes J)). \ar{u}\end{tikzcd}\] From ${\rm Hom}_{{\rm RingRep}_{/(A, G)}}((A, G), (A \oplus J, G \ltimes J))$, we extract the derivation $d$, the crossed homomorphism $\delta$, and the compatibility condition (1) of Definition \ref{def:logder}. Since $M \ltimes J := M \times_{G} (G \ltimes J)$, we have that \[{\rm Hom}_{{\rm LogCrMod}_{/(M \to G)}}(M \to G, M \ltimes J \to G \ltimes J) \cong {\rm Hom}_{{\rm Grp}_{/G}}(G, G \ltimes J).\] In other words, the right-hand ${\rm Hom}$-set is the subset of the set of sections $(A, G) \to (A \oplus J, G \ltimes J)$ in ${\rm RingRep}$ with the property that the two compositions \[M \xrightarrow{\alpha} (A, \cdot) \xrightarrow{({\rm id}_A, d)} (A \oplus J, \cdot)\] and \[M \cong M \times_G G \xrightarrow{{\rm id}_M \times_{{\rm id}_G} ({\rm id}, \delta)} M \times_G (G \ltimes J) =: M \ltimes J \xrightarrow{\alpha \ltimes J} (A \oplus J, \cdot)\] agree. But this is precisely condition (2) of Definition \ref{def:logder}.
\end{proof}

\begin{corollary} For $(A, M)$ a pre-log ring and $J$ an $(A, M)$-bimodule, we have \[{\rm HH}^1((A, M) ; J) \cong {\rm LogDer}((A, M), J) / {\rm InnDer}((A, M), J),\] naturally in $(A, M)$ and $J$. 
\end{corollary}

\begin{proof} The classical proof goes through: The short exact sequence \[0 \to I_{(A, M)} \to (A, M)^e \to A \to 0\] gives rise to the short exact sequence \[0 \to J / Z_A(J)^G \xrightarrow{j \mapsto (d_j, \delta_j)} {\rm LogDer}((A, M), J) \to {\rm Ext}^1_{(A, M)^e}(A, J) \to 0\] upon applying ${\rm Hom}_{(A, M)^e}(-, J)$ by (the proof of) Lemma \ref{lem:loghhzero} and Lemma \ref{lem:unilogder}. The image of the first map is precisely the inner log derivations. 
\end{proof}

\subsection{The exactified abelianization functor} We are now in position to finalize the proof of Theorem \ref{thm:abgrpob}. Let us first assemble the expected description of abelian group objects over a fixed log ring $(A, \varphi \colon M \to G)$ as split square-zero extensions: 

\begin{lemma}\label{lem:abgrpsplit} Let $(p, p^\flat, p^\#) \colon (B, \psi \colon N \to H) \to (A, \varphi \colon M \to G)$ be an abelian group object in ${\rm Log}^{\rm ex}_{/(A, M)}$. There exists an $(A, M)$-bimodule $J$ and an isomorphism \[(A \oplus J, \varphi \ltimes J \colon M \ltimes J \to G \ltimes J) \xrightarrow{\cong} (B, \psi \colon N \to H)\] over and under $(A, \varphi \colon M \to G)$. 
\end{lemma}

\begin{proof} By the definition of abelian group objects, the structure map $(p, p^\flat, p^\#)$ admits a (unit) section. We may thus apply Lemma \ref{lem:retractivelog} to obtain an isomorphism \[(B, \varphi \ltimes {\rm ker}(p^\times) \colon M \ltimes {\rm ker}(p^\times) \to G \ltimes {\rm ker}(p^\times)) \xrightarrow{\cong} (B, \psi \colon N \to H)\] over and under $(A, \varphi \colon M \to G)$. Since $B$ is, in particular, an abelian group object in ${\rm Ring}_{/A}$, there exists an $(A, A)$-bimodule $J$ such that $B \cong A \oplus J$, and $J$ carries the $(A, M)$-bimodule structure specified in the proof of Lemma \ref{lem:abphi}. Now ${\rm ker}(p^\times) = J$, and it is a straight-forward check that this is indeed canonically isomorphic to the split square-zero extension of Construction \ref{constr:splitsqzero}.   
\end{proof}

\begin{proof}[Proof of Theorem \ref{thm:abgrpob}] It remains to describe the left adjoint to the composite \[{\rm Ab}({\rm Log}^{\rm ex}_{/(A, M)}) \to {\rm Log}^{\rm ex}_{/(A, M)} \to {\rm Log}_{/(A, M)}.\] So we study the ${\rm Hom}$-set \[{\rm Hom}_{{\rm Log}_{/(A, M)}}((B, N), (A \oplus J, M \ltimes J))\] for an arbitrary $(A, M)^e$-module $J$; by Lemma \ref{lem:abgrpsplit}, the codomain variable describes an arbitrary object in ${\rm Ab}({\rm Log}^{\rm ex}_{/(A, M)})$. Pulling back along $(B, N) \to (A, M)$ and applying Lemma \ref{lem:logdersqzero}, this ${\rm Hom}$-set is naturally isomorphic to \[{\rm Hom}_{{\rm Log}_{/(B, N)}}((B, N), (B \oplus J, N \ltimes J)) \cong {\rm LogDer}((B, N), J),\] where $J$ is considered a $(B, N)$-bimodule by restriction of scalars along $(B, N) \to (A, M)$. By Lemma \ref{lem:unilogder}, this is naturally isomorphic to \[{\rm Hom}_{{\rm LMod}((B, N)^e)}(I_{(B, N)}, J) \cong {\rm Hom}_{{\rm LMod}((A, M)^e)}((A, M)^e \otimes_{(B, N)^e} I_{(B, N)}, J),\] which concludes the proof. 
\end{proof}

\section{Log Hochschild Cohomology and Square-zero Extensions}\label{sec:lastsection} We now work towards the proofs of Theorems \ref{thm:loghch} and \ref{thm:pullback} from the introduction. Both statements (implicitly) use animated log rings: 

\subsection{Animated pre-log rings} There is a forgetful functor \[{\rm PreLog} \to {\rm Set} \times {\rm Set} \times {\rm Set}, (A, \varphi \colon M \to G) \mapsto (A, M, G).\] This admits a left adjoint ${\rm Free}_{\rm PreLog}(-, -, -)$ by \cite[Remark 14.12(1)]{ARV11}.  

Let ${\rm PolyPreLog}$ denote the category of pre-log rings of the form ${\rm Free}_{\rm PreLog}(S, T, U)$ for finite sets $S, T,$ and $U$. The restricted Yoneda embedding \[y \colon {\rm PreLog} \to {\rm Fun}^{\times}({\rm PolyPreLog}^{\rm op}, {\rm Set})\] to product-preserving presheaves on ${\rm PolyPreLog}^{\rm op}$ is an equivalence of categories.

We define the category of \emph{animated pre-log rings} \[{\rm AniPreLog} := {\rm Fun}^{\times}({\rm PolyPreLog}^{\rm op}, {\rm Ani})\] as the category of product-preserving anima-valued presheaves. By \cite[5.5.9.1, 5.5.9.3]{Lur09}, ${\rm AniPreLog}$ is equivalent to the category simplicial pre-log rings.

\subsection{Animated log rings} We say that an animated pre-log ring $(A, \varphi \colon M \to G)$ is \emph{log} if $\alpha^{-1}A^\times \to A^\times$ is an equivalence, where $A^\times := A \times_{\pi_0(A)} \pi_0(A)^\times$. We follow the strategy of Sagave--Sch\"urg--Vezzosi \cite[Section 3.1]{SSV16} to characterize the log condition in terms of a lifting property. To do so, we shall make use of a map \begin{equation}\label{eq:logliftmap}{\rm Free}_{\rm PreLog}(\{s\}, \emptyset, \emptyset) \to {\rm Free}_{\rm PreLog}(\emptyset, \{t\}, \emptyset).\end{equation} One readily checks that the domain is $({\Bbb Z}[s], {\rm id} \colon \{1\} \to \{1\})$. The map \eqref{eq:logliftmap} is defined to send $s$ to the ring element corresponding to the monoid generator $t$ in the codomain. We shall consider the pushout $P$ of \eqref{eq:logliftmap} along the localization \[{\rm Free}_{{\rm PreLog}}(\{s\}, \emptyset, \emptyset) = ({\Bbb Z}[s], {\rm id} \colon \{1\} \to \{1\}) \to ({\Bbb Z}[s^{\pm 1}], {\rm id} \colon \{1\} \to \{1\}).\] The following is the non-commutative analog of the ideas in \cite[Section 3.1]{SSV16}:

\begin{lemma}\label{lem:contractiblelift} A pre-log ring $(B, \psi \colon N \to H)$ is log if and only if it is local with respect to the map $({\Bbb Z}[s^{\pm 1}], {\rm id} \colon \{1\} \to \{1\}) \to P$; that is, when the space of lifts \[\begin{tikzcd}({\Bbb Z}[s^{\pm 1}], {\rm id} \colon \{1\} \to \{1\}) \ar{d} \ar{r} & (B, \psi \colon N \to H) \\ P \ar[dashed]{ur}\end{tikzcd}\] is contractible.
\end{lemma}

\begin{proof} By construction, the space of maps $P \to (B, \psi \colon N \to H)$ is precisely $\beta^{-1}B^\times$, where $\beta \colon N \to (B, \cdot)$ is the structure map, while the space of maps $({\Bbb Z}[s^{\pm 1}], {\rm id} \colon \{1\} \to \{1\}) \to (B, \psi \colon N \to H)$ is $B^\times$. So having a contractible space of lifts means exactly that $\beta^{-1}B^\times \to B^\times$ is an equivalence. 
\end{proof}

\begin{definition} The category of \emph{animated log rings} ${\rm AniLog}$ is the localization of ${\rm AniPreLog}$ with respect to the map $({\Bbb Z}[s^{\pm 1}], {\rm id} \colon \{1\} \to \{1\}) \to P$. We shall write \[(-)^a \colon {\rm AniPreLog} \to {\rm AniLog}\] for the resulting localization functor. 
\end{definition}

If $(A, \varphi \colon M \to G)$ is an animated pre-log ring, then we obtain a pre-log ring $(\pi_0(A), \pi_0(\varphi) \colon \pi_0(M) \to \pi_0(G))$, and this process preserves the log condition. 

\subsection{The (animated) log enveloping algebra} It is now straight-forward to extend the construction of the log enveloping algebra to the animated context. To distinguish it from the discrete variant of the log enveloping algebra discussed so far, we denote the latter by $(A, M)^e_{\heartsuit}$. 

\begin{definition} We define the \emph{log enveloping algebra} to be the functor \[(-, -)^e \colon {\rm AniPreLog} \to {\rm AniRing}\] obtained by left Kan extending ${\rm PolyPreLog} \to {\rm AniRing}, (P, Q) \mapsto (P, Q)^e_{\heartsuit}$. 
\end{definition}

Since the functor ${\rm PreLog} \to {\rm Ring}, (A, M) \mapsto (A, M)^e_{\heartsuit}$ commutes with sifted colimits, we have that $\pi_0((A, M)^e) \cong (\pi_0(A), \pi_0(M))^e_{\heartsuit}$ for a general animated pre-log ring $(A, M)$. For polynomial pre-log rings $(P, Q)$ resolving the animated pre-log ring $(A, M)$, the $(P, Q)^e_{\heartsuit}$-module maps $(P, Q)^e_{\heartsuit} \to P$ realize to an $(A, M)^e$-module map $(A, M)^e \to A$, and we now write $I_{(A, M)}$ for its fiber. For $J$ a connective left module over $(A, M)^e$, we write $(A, M) \ltimes J$ for the levelwise extension of Construction \ref{constr:splitsqzero} (obtained by modeling ${\rm AniPreLog}$ by simplicial pre-log rings). 

\begin{lemma}\label{lem:unilogderder} Let $(A, M)$ be an animated pre-log ring and let $J$ be a connective left $(A, M)^e$-module. There is a canonical equivalence \[{\rm Map}_{{\rm LMod((A, M)^e)}}(I_{(A, M)}, J) \simeq {\rm Map}_{{\rm AniPreLog_{/(A, M)}}}((A, M), (A, M) \ltimes J)\] of anima. 
\end{lemma}

\begin{proof} This is analogous to \cite[Lemma 4.6]{SSV16} in the classically commutative setting: It follows from the discrete statement, which is Lemmas \ref{lem:unilogder} and \ref{lem:logdersqzero}. 
\end{proof}

\subsection{Square-zero extensions and derivations} Recall that, if $B \to A$ is a square-zero extension of rings with kernel $J$, then $B \times_A B \cong B \oplus J$. This fact has the following logarithmic analog:

\begin{lemma}\label{lem:sqzeropull} Let $(p, p^\flat, p^\#) \colon (B, \psi \colon N \to H) \to (A, \varphi \colon M \to G)$ be an exact square-zero extension of log rings by the $(A, M)^e_{\heartsuit}$-module $J$. The map \[(B \oplus J, \psi \ltimes {\rm id} \colon N \ltimes J \to H \ltimes J) \to (B, \psi \colon N \to H) \times_{(A, \varphi \colon M \to G)} (B, \psi \colon N \to H)\] given by \[(b, j) \mapsto (b, b + j), \quad (n, j) \mapsto (n, \beta^{-1}(1 + j)n), \quad (h, j) \mapsto (h, \psi(\beta^{-1}(1 + j))h)\] at the level of rings, monoids, and groups, is an isomorphism. 
\end{lemma}

\begin{proof}  The projection to the first factor of the pullback is exact, and it is split by the diagonal. Its kernel on units consists of pairs $(1, 1 + j)$, which we identify with ${\rm ker}(p_{|B^\times})$. By inspection, Lemma \ref{lem:retractivelog} applies to conclude. 
\end{proof}

\begin{proof}[Proof of Theorem \ref{thm:pullback}] The following argument is similar to that of Pridham \cite[Lemma 1.8]{Pri24} in the context of simplicial commutative rings. Consider the square \begin{equation}\label{eq:thmpullbacksquare}\begin{tikzcd}(B, N) \ar{r}{\Delta} \ar{d} & (B, N)^{\times_{(A, M)} \bullet + 1} \ar{d} \\ (A, M) \ar{r}{(d_0, d_0^\flat, d_0^\#)} & (A, M) \ltimes J^{\oplus \bullet}\end{tikzcd}\end{equation} of simplicial pre-log rings, where $(B, N)^{\times_{(A, M)} \bullet + 1}$ is the \v{C}ech nerve of the given map $(B, N) \to (A, M)$, while $(A, M) \ltimes J^{\oplus \bullet}$ is the bar construction of $(A, M) \ltimes J$ in ${\rm Ab}({\rm Log}^{\rm ex}_{/(A, M)})$; that is, it is $(A, M) \ltimes B_{\bullet}J$. Under the degreewise isomorphisms $(B, N)^{\times_{(A, M)} \bullet + 1} \cong (B, N) \ltimes J^{\oplus \bullet}$ provided by Lemma \ref{lem:sqzeropull}, the right-hand vertical map in \eqref{eq:thmpullbacksquare} is simply induced by $(B, N) \to (A, M)$. 

The square \eqref{eq:thmpullbacksquare} is equivalent to the square \begin{equation}\label{eq:thmpullbacksquare2}\begin{tikzcd}(B, N) \ar{r} \ar{d} & (A, M) \ar{d}{(d, d^\flat, d^\#)} \\ (A, M) \ar{r}{(d_0, d_0^\flat, d_0^\#)} & (A, M) \ltimes J[1]\end{tikzcd}\end{equation} which we claim to be cartesian. For this, we apply the Bousfield--Friedlander theorem \cite[Theorem B.4]{BF78}. Here, \eqref{eq:thmpullbacksquare} is considered one of bisimplicial objects that are constant in the vertical direction, rendering the $\pi_*$-Kan condition vacuous. The square \eqref{eq:thmpullbacksquare} is levelwise cartesian, and the right-hand vertical map of \eqref{eq:thmpullbacksquare} is a Kan fibration on vertical path components: On the ring and group components this follows from surjectivity of $B \to A$ and $H \to G$, while on the monoid component this follows from exactness. Thus \cite[Theorem B.4]{BF78} indeed applies to infer that \eqref{eq:thmpullbacksquare2} is cartesian in the ring, monoid, and group component separately, rendering \eqref{eq:thmpullbacksquare2} cartesian. 

Conversely, suppose that the derivation $(d, d^\flat, d^\#) \colon (A, M) \to (A, M) \ltimes J[1]$ is given and that $(B, N)$ is defined to sit in the pullback diagram \eqref{eq:thmpullbacksquare2}. It is then the case that $B \to A$ is a square-zero extension by $J$. To argue that the map $(B, N) \to (A, M)$ is strict, we first observe that it is exact since $(d, d^\flat, d^\#)$ is, so that Lemma \ref{lem:exstrict} applies. In this case, \eqref{eq:thmpullbacksquare2} is equivalent to the square \[\begin{tikzcd}(B, N) \ar{r}{\Delta} \ar{d} & (B, N)^{\times_{(A, M)} \bullet + 1} \ar{d} \\ (A, M) \ar{r}{} & (A, M)^{\times_{(A, M) \ltimes J[1]} \bullet + 1},\end{tikzcd}\] where the lower right-hand corner is the \v{C}ech nerve associated to the zero section $(d_0, d_0^\flat, d_0^\#) \colon (A, M) \to (A, M) \ltimes J[1]$ and the right-hand vertical map is induced by $(B, N) \to (A, M)$ and $(d, d^\flat, d^\#)$. Under the isomorphism $(B, N)^{\times_{(A, M)} \bullet + 1} \cong (B, N) \ltimes J^{\oplus \bullet}$ provided by Lemma \ref{lem:sqzeropull} and the equivalence $(A, M)^{\times_{(A, M) \ltimes J[1]} \bullet + 1} \simeq (A, M) \ltimes B_{\bullet}J$, the right-hand vertical map corresponds to that induced by $(B, N) \to (A, M)$. Thus, any derivation arises from the procedure \eqref{eq:thmpullbacksquare}, as desired. 
\end{proof}

\begin{proof}[Proof of Theorem \ref{thm:loghch}] By Theorem \ref{thm:pullback}, the isomorphism classes of the relevant square-zero extensions identify with \[\pi_0{\rm Map}_{{\rm AniLog}_{/(A, M)}}((A, M), (A, M) \ltimes J[1]) \cong \pi_0{\rm Map}_{{\rm LMod}((A, M)^e)}(I_{(A, M)}, J[1]),\] where the isomorphism is due to Lemma \ref{lem:unilogderder}. Combining the fact that $J$ is discrete with the long exact sequence associated to the defining fiber sequence $I_{(A, M)} \to (A, M)^e \to A$ identifies this further with \[\pi_0{{\rm Map}_{{\rm LMod}((A, M)^e)}}(A, J[2]) \cong {\rm Ext}^2_{(A, M)^e}(A, J),\] as desired. 
\end{proof}

\vspace{10 mm}

\end{document}